\documentclass[11pt,a4paper]{article}
\usepackage{amsmath,amsfonts,amssymb,amsthm}
\usepackage{dsfont,mathrsfs,stmaryrd,enumerate}

\usepackage{curves,epic,eepic}
\usepackage{graphicx}

\newtheorem{defi}{Definition}[section]
\newtheorem{lemma}[defi]{Lemma}
\newtheorem{prop}[defi]{Proposition}
\newtheorem{propdefi}[defi]{Proposition--Definition}
\newtheorem{thm}[defi]{Theorem}
\newtheorem{cor}[defi]{Corollary}
\theoremstyle{definition}
\newtheorem{rem}[defi]{Remark}

\def\Z{\mathds Z}
\def\Q{\mathds Q}

\def\C{\mathds C}
\def\e{\varepsilon}
\def\phi{\varphi} 
\def\rho{\varrho}
\def\sub{\subseteq}
\def\id{{\it id}}
\def\A{\mathds A}
\def\P{\mathds P}
\def\O{\mathcal O}
\def\<{\left<}
\def\>{\right>}

\DeclareMathOperator{\Hom}{Hom}
\DeclareMathOperator{\Aut}{Aut}
\DeclareMathOperator{\End}{End}
\DeclareMathOperator{\SHom}{{\mathscr H}\!\!\textit{om}}

\newcommand{\catXschemes}[1]{(\textit{$#1$-schemes})}
\def\catsets{(\textit{sets})}

\DeclareMathOperator{\GL}{GL}
\DeclareMathOperator{\SL}{SL}
\DeclareMathOperator{\Spec}{Spec}
\DeclareMathOperator{\rk}{rk}
\DeclareMathOperator{\ord}{ord}

\DeclareMathOperator{\KG}{\textit{KG}}
\DeclareMathOperator{\LG}{\textit{LG}}

\DeclareMathOperator{\BD}{\textit{BD}}
\DeclareMathOperator{\BT}{\textit{BT}}
\DeclareMathOperator{\BO}{\textit{BO}}
\DeclareMathOperator{\BI}{\textit{BI}}

\DeclareMathOperator{\Hilb}{Hilb}

\DeclareMathOperator{\GHilb}{G-Hilb}
\DeclareMathOperator{\FGHilb}{\underline{G-Hilb}}
\newcommand{\XHilb}[1]{{#1}{\textrm{-}}\Hilb}

\DeclareMathOperator{\Irr}{Irr}
\DeclareMathOperator{\graph}{Graph}
\newcommand{\HSymb}[2]{(\!(#1,#2)\!)}

\begin{document}

\title{\vspace{-9mm}
McKay correspondence over non algebraically
closed fields}
\author{Mark Blume}
\date{\vspace{-6mm}}

\maketitle

\begin{abstract}
The classical McKay correspondence for finite subgroups 
$G$ of $\SL(2,\C)$ gives a bijection between isomorphism classes 
of nontrivial irreducible representations of $G$ and irreducible
components of the exceptional divisor in the minimal resolution 
of the quotient singularity $\A^2_\C/G$. 
Over non algebraically closed fields $K$ there may exist
representations irreducible over $K$ which split over $\overline{K}$.
The same is true for irreducible components of the exceptional
divisor. In this paper we show that these two phenomena
are related and that there is a bijection between nontrivial
irreducible representations and irreducible components of the 
exceptional divisor over non algebraically closed fields $K$ 
of characteristic $0$ as well.
\end{abstract}

\section{Introduction}

Let $G$ be a finite group operating on a smooth variety $M$ 
over $\C$, e.g.~$M=\A^n_\C$ and a linear operation of a 
finite subgroup $G\subset\SL(n,\C)$.
Usually the quotient $M/G$ is singular and one considers  
resolutions of singularities $Y\to M/G$ with some minimality 
property. 
A method to construct resolutions of quotient singularities
is the $G$-Hilbert scheme $\GHilb M$ introduced in \cite{ItNm96}, 
\cite{ItNm99}. Under some conditions the $G$-Hilbert scheme is 
irreducible, nonsingular and $\GHilb M\to M/G$ a crepant resolution 
\cite{BKR01}. In particular, this applies to the operation of 
finite subgroups $G\subset\SL(n,\C)$ on $\A^n_\C$ for $n\leq 3$.
For $G\subset\SL(2,\C)$ there are also other methods to show that 
the $G$-Hilbert scheme is the minimal resolution, see \cite{ItNm96},
\cite{ItNm99}.\medskip 

The McKay correspondence in general describes the resolution 
$Y$ in terms of the representation theory of the group $G$, 
see \cite{Re97}, \cite{Re99} for expositions of this subject. 
Part of the correspondence for $G\subset\SL(2,\C)$ is a bijection 
between irreducible components of the exceptional divisor $E$ 
and isomorphism classes of nontrivial irreducible representations 
of the group $G$ and moreover an isomorphism of graphs between 
the intersection graph of components of $E$ and the 
representation graph of $G$, both being graphs of ADE type. 
This was the observation of McKay \cite{McK80}.\medskip

The new contribution in this paper is to consider McKay correspondence
over non algebraically closed fields. We will work over a field $K$ 
that is not assumed to be algebraically closed but always of 
characteristic $0$ and extend the McKay correspondence to this 
slightly more general situation.
Over non algebraically closed $K$ it is natural to consider finite 
group schemes instead of simply finite groups.
In comparison with the situation over algebraically closed fields
there may exist both representations of $G$ and components of $E$ 
that are irreducible over $K$ but split over the algebraic closure. 
We will see that these two kinds of splitting that arise by 
extending the ground field are related by investigating the
operation of the Galois group. For this we introduce Galois-conjugate 
representations and consider the Galois operation on the 
$G$-Hilbert scheme.
The following McKay correspondence over arbitrary fields $K$ 
of characteristic $0$ will be consequence of more detailed 
theorems in section 5.

\begin{thm} Let $K$ be any field of characteristic $0$ and
$G\subset\SL(2,K)$ a finite subgroup scheme. Then there is a 
bijection between the set of irreducible components of the 
exceptional divisor $E$ and the set of isomorphism classes 
of nontrivial irreducible representations of $G$ and moreover 
an isomorphism between the intersection graph of the irreducible 
components of $E_{\rm red}$ and the representation graph of $G$.
\end{thm}

Examples are discussed in subsection \ref{subsec:ex}, the possible 
graphs can be found in subsection \ref{subsec:repgraphs}.
As already observed in \cite{Li69}, considering the rational double 
points over non algebraically closed fields one finds the remaining 
Dynkin diagrams of types $(B_n)$, $(C_n)$, $(F_4)$, $(G_2)$.
The methods of this paper should also apply to other situations,
in particular to the McKay correspondence for finite small subgroups 
of $\GL(2,\C)$ and give a similar generalisation as in the $\SL$-case.

\medskip

This paper is organised as follows.
Section 2 shortly summarises some techniques used in this paper,
namely $G$-sheaves for group schemes $G$ and $G$-Hilbert 
schemes.
Section 3 is concerned with the relations between Galois operations 
and decompositions into irreducible components both of schemes and 
representations. We introduce the notion of Galois-conjugate 
representations and $G$-sheaves and we describe the Galois operation 
on $G$-Hilbert schemes.  
In section 4 we collect some data of the finite subgroup schemes of 
$\SL(2,K)$ and list possible representation graphs. In addition we 
investigate under what conditions a finite subgroup of $\SL(2,C)$, 
$C$ the algebraic closure of $K$, is realisable as a subgroup 
of $\SL(2,K)$.
Section 5 contains the theorems of McKay correspondence over
non algebraically closed fields. We consider two constructions,
the stratification of the $G$-Hilbert scheme and the tautological 
sheaves, originating from \cite{ItNm96} and \cite{GV83} 
respectively, that are known to give a McKay correspondence 
over $\C$ and formulate them for not necessarily algebraically 
closed $K$.\medskip

{\it Acknowledgements.} The suggestion to investigate McKay 
correspondence over non algebraically closed fields is due 
to Victor Batyrev.\medskip

{\it Notations.} In general we write a lower index for base 
extensions, for example if $X,T$ are $S$-schemes then $X_T$ denotes 
the $T$-scheme $X\times_ST$ or if $V$ is a representation over a 
field $K$ then $V_L$ denotes the representation $V\otimes_KL$ over 
the extension field $L$. Likewise, if $\phi\colon X\to Y$ is a morphism 
of $S$-schemes, we write $\phi_T\colon X_T\to Y_T$ for its base extension
with respect to $T\to S$.

\section{Preliminaries}

\subsection{G-sheaves}

Let $K$ be a field. Let $G$ be a group scheme over $K$ 
with $p\colon G\to\Spec K$ the projection, $e\colon\Spec K\to G$ 
the unit, and $m\colon G\times_KG\to G$ the multiplication. 
For affine $G=\Spec A$, $A$ has the structure of a Hopf 
algebra over $K$, the coalgebra structure being equivalent 
to the group structure of $G$.\medskip

Let $X$ be a $G$-scheme over $K$, that is a $K$-scheme with
an operation $s_X\colon G\times_KX\to X$ of the group scheme $G$
over $K$.
We have to use a more general notion of a $G$-sheaf than in  
\cite{BKR01}, we adopt the definition of \cite{Mu}:
a (quasicoherent, coherent) $G$-sheaf on $X$ is a 
(quasicoherent, coherent) $\O_X$-module $\mathscr F$ with 
an isomorphism $\lambda^{\mathscr F}\colon s_X^*\mathscr F
\stackrel{\!_\sim}{\longrightarrow}p_X^*\mathscr F$ of $\O_{G\times_K X}$-modules
satisfying the conditions
(i) the restriction of $\lambda^{\mathscr F}$ to the unit in $G_X$
is the identity, i.e. $e_X^*\lambda^{\mathscr F}\colon e_X^*s_X^*
\mathscr F\to e_X^*p_X^*\mathscr F$ identifies with
$\id_{\mathscr F}\colon\mathscr F\to\mathscr F$,
and (ii) $(m\times\id_X)^*\lambda^{\mathscr F}
=p_{23}^*\lambda^{\mathscr F}\circ(\id_G\times s_X)^*
\lambda^{\mathscr F}$, where $p_{23}\colon G\times_K G\times_K 
X\to G\times_K X$ is the projection to the factors $2$ and $3$.

\begin{rem}\label{rem:Gsheaves} 
We summarise relevant properties of $G$-sheaves.\\
(1) There is the canonical notion of $G$-equivariant homomorphisms
between $G$-sheaves $\mathscr F,\mathscr G$ on $X$, 
the set of these is denoted by $\Hom^G_X(\mathscr F,\mathscr G)$.
Kernels and cokernels of $G$-equivariant homomorphisms have natural 
$G$-sheaf structures.\\
(2) Assume $G=\Spec A$ affine and let $X$ be a $G$-scheme with 
trivial $G$-operation, i.e. $s_X=p_X$. 
Then the $G$-sheaf structure of a $G$-sheaf $\mathscr F$ is equivalent 
to a homomorphism of $\O_X$-modules $\rho\colon\mathscr F\to A\otimes_K
\mathscr F$ satisfying the usual conditions of a comodule. This 
relation can be constructed using the adjunction $({p_X}^*,p_{X*})$.
Further, notions such as ``subcomodule'', ``homomorphism of comodules'', 
etc.\ correspond to ``$G$-subsheaf'', ``equivariant homomorphism'', etc.. 
The $G$-invariant part $\mathscr F^G\sub\mathscr F$ is defined by 
$\mathscr F^G(U):=\{f\in\mathscr F(U)\:|\:\rho(f)=1\otimes f\}$ 
for open $U\sub X$.\\
(3) For an $A$-comodule $\mathscr F$ on $X$ a decomposition 
of $A$ into a direct sum $A=\bigoplus_i A_i$ of subcoalgebras 
$A_i$ determines a direct sum decomposition 
$\mathscr F=\bigoplus_i\mathscr F_i$ into subcomodules 
(take preimages $\rho^{-1}(A_i\otimes_K\mathscr F)$), where 
the comodule structure of $\mathscr F_i$ reduces to an 
$A_i$-comodule structure.\\
(4) A $G$-sheaf on $X=\Spec K$ (or an extension field of $K$) 
we also call a representation. Dualisation of an $A$-comodule 
$V$ over $K$ leads to a $\KG$-module $V^\vee$, where 
$\KG=A^\vee=\Hom_K(A,K)$ with algebra structure dual to the 
coalgebra structure of $A$.\\
(5) For quasicoherent $G$-sheaves $\mathscr F,\mathscr G$ with
$\mathscr F$ finitely presented the sheaf 
$\SHom_{\O_X}(\mathscr F,\mathscr G)$ carries a natural $G$-sheaf 
structure. For locally free $\mathscr F$ one defines the dual 
$G$-sheaf by $\mathscr F^\vee=\SHom_{\O_X}(\mathscr F,\O_X)$. 
In the case of trivial $G$-operation on $X$ there is the 
component $\SHom^G_{\O_X}(\mathscr F,\mathscr G)$ of 
$\SHom_{\O_X}(\mathscr F,\mathscr G)$, the sheaf of 
equivariant homomorphisms, that can either be described as
$G$-invariant part $(\SHom_{\O_X}(\mathscr F,\mathscr G))^G\sub
\SHom_{\O_X}(\mathscr F,\mathscr G)$ or by 
$\SHom^G_{\O_X}(\mathscr F,\mathscr G)(U)=\Hom^G_U(\mathscr F|_U,
\mathscr G|_U)$ for open $U\sub X$.\\
(6) Functors for sheaves like $\otimes, f^*, \ldots$ as well
have analogues for $G$-sheaves, e.g.~for equivariant $f\colon Y\to X$ 
and a $G$-sheaf $\mathscr F$ on $X$ the sheaf $f^*\mathscr F$ 
has a natural $G$-sheaf structure.\\
(7) Natural isomorphisms for sheaves lead to isomorphisms 
for $G$-sheaves, e.g.\ under some conditions there is an 
isomorphism $f^*\SHom_{\O_X}(\mathscr F,\mathscr G)
\cong\SHom_{\O_Y}(f^*\mathscr F,f^*\mathscr G)$ 
and this isomorphism becomes an isomorphism of $G$-sheaves
provided that $f$ is equivariant and $\mathscr F,\mathscr G$
are $G$-sheaves.
Other examples are $f^*(\mathscr F\otimes_{\O_X}\mathscr G)
\cong f^*\mathscr F\otimes_{\O_Y}f^*\mathscr G$ and
$\SHom_{\O_X}(\mathscr F\otimes_{\O_X}\mathscr E,
\mathscr G)\cong\SHom_{\O_X}(\mathscr F,\mathscr E^\vee
\otimes_{\O_X}\mathscr G)$.\\
(8) Base extension $K\to L$ makes out of a $G$-scheme $X$ over $K$
a scheme $X_L$ with a $G$-scheme or a $G_L$-scheme structure, 
the operation given by $s_{X_L}=(s_X)_L$.
A $G$-sheaf $\mathscr F$ on a $G$-scheme $X$ gives rise to 
a $G$-sheaf $\mathscr F_L=\mathscr F\otimes_KL=f^*\mathscr F$
on $X_L$, where $f\colon X_L\to X$. $\mathscr F_L$ can be considered 
as a $G_L$-sheaf on the $G_L$-scheme $X_L$ over $L$.
\end{rem}

\subsection{G-Hilbert schemes}

Let $G=\Spec A$ be a finite group scheme over a field $K$, 
assume that its Hopf algebra $A$ is cosemisimple (that is, $A$ is 
sum of its simple subcoalgebras, see \cite[Ch. XIV]{Sw} and 
subsection \ref{subsec:cosemisimple-galoisext} below).\medskip

For us the $G$-Hilbert scheme $\GHilb_KX$ of a $G$-scheme $X$ over 
$K$ will be by definition the moduli space of $G$-clusters,
i.e. parametrising $G$-stable finite closed subschemes $Z\sub X_L$, 
$L$ an extension field of $K$, with $H^0(Z,\O_Z)$ isomorphic to the 
regular representation of $G$ over $L$. 
We recall its construction (a variation of the Quot scheme construction
of \cite{Gr61}), for a detailed discussion including the generalisation 
to finite group schemes with cosemisimple Hopf algebra over arbitrary 
base fields see \cite{Bl11}.\medskip

Let $X$ be a $G$-scheme algebraic over $K$, assume that a geometric 
quotient $\pi\colon X\to X/G$, $\pi$ affine, exists. Then the $G$-Hilbert 
functor $\FGHilb_KX\colon\catXschemes{K}^\circ\to\catsets$ given by
\[
\FGHilb_KX(T):=
\left\{
\begin{array}{l}
\textit{Quotient $G$-sheaves}\;\;[0\to\mathscr I\to\O_{X_T}\to
\O_Z\to 0]\;\;\textit{on $X_T$},\\
Z\;\textit{finite flat over $T$},\;\;
\textit{for $t\in T$: $H^0(Z_t,\O_{Z_t})$ isomorphic}\\
\textit{to the regular representation}\\
\end{array}
\right\}
\]
is represented by an algebraic $K$-scheme $\GHilb_KX$. Here we 
write $[0\to\mathscr I\to\O_{X_T}\to\O_Z\to 0]$ for an exact 
sequence $0\to\mathscr I\to\O_{X_T}\to\O_Z\to 0$ of quasicoherent 
$G$-sheaves on $X_T$ with $\mathscr I,\O_Z$ specified up to 
isomorphism, that is either a quasicoherent $G$-subsheaf 
$\mathscr I\sub\O_{X_T}$ or an equivalence class $[\O_{X_T}\to\O_Z]$ 
of surjective equivariant homomorphisms of quasicoherent $G$-sheaves 
with two of them equivalent if their kernels coincide.\\
There is the natural morphism $\tau\colon\GHilb_KX\to X/G$, which is 
projective and as a map of points takes $G$-clusters to the 
corresponding orbits.\medskip 

In this paper we are interested in the case $G\subset\SL(2,K)$ 
operating on $X=\A_K^2$ over fields $K$ of characteristic $0$. 

\begin{prop}
The $G$-Hilbert scheme $\GHilb_K\A^2_K$ is irreducible and 
nonsingular. The morphism $\tau\colon\GHilb_K\A^2_K\to\A^2_K/G$ 
is birational and the minimal resolution of $\A^2_K/G$.
\end{prop}
\begin{proof}
This is known for algebraically closed fields of characteristic 
$0$ \cite{ItNm96}, \cite{ItNm99}, \cite{BKR01}.
From this the statements about irreducibility and nonsingularity 
for not necessarily algebraically closed $K$ follow, use that
for $C$ the algebraic closure $(\GHilb_K\A^2_K)_C\cong
\XHilb{G_C}_C\A^2_C$ (see \cite{Bl11}).
The morphism $\tau\colon\GHilb_K\A^2_K\to\A^2_K/G$ is known to be 
birational.
The base extension $(\GHilb_K\A^2_K)_C\to(\A^2_K/G)_C$ 
identifies with the natural morphism $\XHilb{G_C}_C\A^2_C
\to\A^2_C/G_C$ (follows directly from the functorial definition 
of $\tau$, see e.g.~\cite{Bl11}).
So the statement about minimality as well follows from the same 
statement for algebraically closed fields.
\end{proof}

\section{Galois operation and irreducibility}

\subsection{(Co)semisimple (co)algebras and Galois extensions}
\label{subsec:cosemisimple-galoisext}

Let $K$ be a field and $K\to L$ a Galois extension, 
$\Gamma:=\Aut_K(L)$.
As reference for simple and semisimple algebras we use 
\cite[Alg\'ebre, Ch. VIII]{B}, for coalgebras and comodules
\cite{Sw}. 
Note that for a $K$-vector space $V$ (maybe with some additional 
structure) $\Gamma$ operates on the base extension $V_L=V\otimes_KL$
via the second factor.

\begin{prop}\label{prop:simplealgKL}
Let $F$ be a simple $K$-algebra. 
Assume that $F_L$ is semisimple, let $F_L=\bigoplus_{i=1}^rF_{L,i}$
be its decomposition into simple components.
Then $\Gamma$ permutes the simple summands $F_{L,i}$ and
the operation on the set $\{F_{L,1},\ldots,F_{L,r}\}$
is transitive.
\end{prop}
\begin{proof}
The $F_{L,i}$ are the minimal two-sided ideals of $F_L$. Since
any $\gamma\in\Gamma$ is an automorphism of $F_L$ as a $K$-algebra
or ring, the $F_{L,i}$ are permuted by $\Gamma$.\\
Let $U=\sum_{\gamma\in\Gamma}\gamma F_{L,1}$ and $V$ the sum
over the remaining $F_{L,i}$. Then $F_L=U\oplus V$, $U$ and $V$ 
are $\Gamma$-stable and thus $U=U'_L$, $V=V'_L$ for $K$-subspaces
$U',V'\sub F$ by \cite[Algebra II, Ch. V, \S 10.4]{B}, since
$K\to L$ is a Galois extension.
It follows that $F=U'\oplus V'$ with $U',V'$ two-sided ideals of $F$.
Since $F$ is simple, $V'=0$, $U=F_L$ and the operation is transitive.
\end{proof}

A coalgebra $C\neq 0$ is called simple, if it has no subcoalgebras 
except $\{0\}$ and $C$. A coalgebra is called cosemisimple, if it 
is the sum of its simple subcoalgebras, in which case this sum is 
direct. For cosemisimple $C$ the simple subcoalgebras are the 
isotypic components of $C$ as a $C$-comodule (left or right), 
so they correspond to the isomorphism classes of simple 
representations of $G$ over $K$. 

\begin{prop}\label{prop:cosemisimple-baseext}
Let $C$ be a finite dimensional coalgebra over $K$. Then $C$ is 
cosemisimple if and only if $C_L$ is cosemisimple.
\end{prop}
\begin{proof}
This is equivalent to the dual statement for finite dimensional 
semisimple $K$-algebras 
\cite[Alg\'ebre, Ch. VIII, \S 7.6, Thm. 3, Cor. 4]{B}.
\end{proof}

For simple coalgebras there is a result similar to proposition 
\ref{prop:simplealgKL} and proven analogously, note that simple 
coalgebras are finite dimensional.

\begin{prop}\label{prop:simplecoalgKL}
Let $C$ be a simple coalgebra over $K$. 
Then $C_L$ is cosemisimple, and if $C_L=\bigoplus_iC_{L,i}$
is its decomposition into simple components, then $\Gamma$ 
transitively permutes the simple summands $C_{L,i}$.\qed
\end{prop}

\begin{cor}\label{cor:decompcoalg} 
Let $C$ be a cosemisimple coalgebra over $K$.
Then $C_L$ is cosemisimple, and if $C=\bigoplus_jC_j$ resp.\ 
$C_L=\bigoplus_iC_{L,i}$ are the decompositions into simple 
subcoalgebras, then:

\vspace{-2mm}

\begin{enumerate}[\rm (i)]\setlength{\itemsep}{0mm}
\item The decomposition  $C_L\cong\bigoplus_iC_{L,i}$ is a 
refinement of the decomposition $C_L\cong\bigoplus_j(C_j)_L$.
\item $\Gamma$ transitively permutes the summands 
$C_{L,i}$ of $(C_j)_L$ for any $j$.\\
Therefore $(C_j)_L=\sum_{\gamma\in\Gamma}\gamma C_{L,i}$,
if $C_{L,i}$ is a summand of $(C_j)_L$.\qed
\end{enumerate}
\end{cor}

This applies to the situation considered in this paper.
Assume that the field $K$ is of characteristic $0$ and let 
$G=\Spec A$ be a finite group scheme over $K$, $|G|:=\dim_KA$. 
Define $\KG$ to be the $K$-vector space $A^\vee=\Hom_K(A,K)$ 
with algebra structure dual to the coalgebra structure of $A$.
In this situation the algebra $A$ is always reduced and for a
suitable algebraic extension field $L$ of $K$ the group scheme 
$G_L$ is discrete. Then $G(L)$ is a finite group of order $|G|$
and the algebra $\LG=(\KG)_L$ is isomorphic to the group algebra 
of the group $G(L)$ over $L$.
By semisimplicity of group algebras for finite groups over fields of 
characteristic $0$ and proposition \ref{prop:cosemisimple-baseext} 
one obtains:

\begin{prop}
Let $G=\Spec A$ be a finite group scheme over a field $K$ of 
characteristic $0$. Then the Hopf algebra $A$ is cosemisimple
and so are its base extensions $A_L$ with respect to field
extensions $K\to L$.\qed
\end{prop}

\subsection{Irreducible components of schemes and Galois extensions}

Let $X$ be a $K$-scheme. For an extension field $L$ of $K$ the group
$\Gamma=\Aut_K(L)$ operates on $X_L=X\times_K\Spec L$ by 
automorphisms of $K$-schemes via the second factor.
For simplicity we denote the morphisms $\Spec L\to\Spec L$, 
$X_L\to X_L$ coming from $\gamma\colon L\to L$ by $\gamma$ as well.\medskip

A point of $X$ may decompose over $L$, this way a point 
$x\in X$ corresponds to a set of points of $X_L$, the 
preimage of $x$ with respect to the projection $X_L\to X$.
In particular this applies to closed points and to irreducible
components. These sets are known to be exactly the $\Gamma$-orbits.

\begin{prop}\label{prop:irredcomp}
Let $X$ be an algebraic $K$-scheme and $K\to L$ be a 
Galois extension, $\Gamma:=\Aut_K(L)$. 
Then points of $X$ correspond to $\Gamma$-orbits of points of $X_L$,
the $\Gamma$-orbits are finite.
\end{prop}
\begin{proof}
Taking fibers, the proposition reduces to the following statement:\smallskip

{\it Let $F$ be the quotient field of a commutative integral 
$K$-algebra of finite type. 
Then $F_L=F\otimes_KL$ has only finitely many prime ideals 
and they are $\Gamma$-conjugate.}\smallskip

Proof. $F_L$ is integral over $F$ because this property is
stable under base extension
\cite[Commutative Algebra, Ch. V, \S 1.1, Prop. 5]{B}. 
It is clear that every prime ideal of $F_L$ lies above the prime 
ideal $(0)$ of $F$. There are no inclusions between the prime 
ideals of $F_L$ \cite[Commutative Algebra, Ch. V, \S 2.1, 
Proposition 1, Corollary 1]{B}.
Since every prime ideal of $F_L$ is a maximal ideal and $F_L$ 
is noetherian (a localisation of an $L$-algebra of finite type), 
$F_L$ is artinian, it has only finitely many prime ideals
$Q_1,\ldots,Q_r$.

$F_L$ has trivial radical \cite[Alg\'ebre, Ch. VIII, \S 7.3, Thm. 1, 
also \S7.5 and \S 7.6, Cor. 3]{B}.
Being an artinian ring without radical, i.e. semisimple
\cite[Alg\'ebre, Ch. VIII, \S 6.4, Thm. 4, Cor. 2 and Prop. 9]{B},
$F_L$ decomposes as a $L$-algebra into a direct sum
\[\textstyle F_L\cong\bigoplus_{i=1}^rF_{L,i}\]
of fields $F_{L,i}\cong F_L/Q_i$ (this can easily be seen directly,
however, it is part of the general theory of semisimple algebras 
developed in \cite[Alg\'ebre, Ch. VIII]{B} that contains the 
representation theory of finite groups schemes with cosemisimple Hopf
algebra as another special case). 

$\Gamma$ operates on $F_L$, it permutes the $Q_i$ and the simple 
components $F_{L,i}$ of $F_L$ transitively by proposition 
\ref{prop:simplealgKL}.
\end{proof}

\subsection{Galois operation on G-Hilbert schemes}

Let $Y$ be a $K$-scheme, $L$ an extension field of $K$
and $\Gamma=\Aut_K(L)$.\medskip
 
For an $L$-scheme $f\colon T\to\Spec L$ and $\gamma\in\Gamma$
define the $L$-scheme $\gamma_* T$ to be the scheme $T$
with structure morphism $\gamma\circ f$.
For a morphism $\alpha\colon T'\to T$ of $L$-schemes let 
$\gamma_*\alpha$ be the same morphism $\alpha$ considered 
as an $L$-morphism $\gamma_*T'\to\gamma_*T$.\\
For a morphism $\alpha\colon Y_L\to Y'_L$ of $L$-schemes and
$\gamma\in\Gamma$ define the conjugate morphism $\alpha^\gamma$ by 
$\alpha^\gamma:=\gamma\circ(\gamma_*\alpha)\circ\gamma^{-1}$, 
which again is a morphism of $L$-schemes. Here 
$\gamma\colon\gamma_*Y_L\to Y_L$ is a morphism over $L$.\medskip

Let $T$ be an $L$-scheme defined over $K$, that is $T=T'_L$ for 
some $K$-scheme $T'$. The group $\Gamma$ operates on the set 
$Y_L(T)$ of morphisms $T\to Y_L$ over $L$ by
\[
\begin{array}{rrcl}
\gamma:&Y_L(T)&\;\to\;&Y_L(T)\\
&\alpha&\;\mapsto\;&\alpha^\gamma=\gamma\circ(\gamma_*\alpha)
\circ\gamma^{-1}\\
\end{array}
\]
Consider the case of $G$-Hilbert schemes. Let $G$ be a finite 
group scheme over $K$, $X$ be a $G$-scheme over $K$ and assume 
that the $G$-Hilbert functor is represented by a $K$-scheme 
$\GHilb_KX$. There is the canonical isomorphism of $L$-schemes
$(\GHilb_KX)_L\cong\XHilb{G_L}_LX_L$ 
(see \cite{Bl11}), obtained by identifying 
$X\times_KT=X_L\times_LT$ for $L$-schemes $T$.

\begin{prop}
Let $T$ be an $L$-scheme defined over $K$.
Then, for a morphism $\alpha\colon T\to\XHilb{G_L}_LX_L$ of $L$-schemes 
corresponding to a quotient $[0\to\mathscr I\to\O_{X_T}\to
\O_Z\to 0]$ and for $\gamma\in\Gamma$, the $\gamma$-conjugate 
morphism $\alpha^\gamma$ corresponds to the quotient 
$[0\to\gamma_*\mathscr I\to\O_{X_T}\to\O_{\gamma Z}\to 0]$.
\end{prop}
\begin{proof} For a morphism of $L$-schemes
$\alpha\colon T\to\XHilb{G_L}_LX_L\cong(\GHilb_KX)\times_K\Spec L$ 
consider the commutative diagram of $L$-morphisms

\medskip

\noindent
\begin{picture}(160,25)
\put(40,20){\makebox(0,0)[c]{$\gamma_*T$}}
\put(100,20){\makebox(0,0)[c]{$(\GHilb_KX)\times_K(\gamma_*\Spec L)$}}
\put(40,5){\makebox(0,0)[c]{$T$}}
\put(96,5){\makebox(0,0)[c]{$(\GHilb_KX)\times_K\Spec L$}}
\put(48,20){\vector(1,0){24}}
\put(48,5){\vector(1,0){24}}
\put(40,17){\vector(0,-1){9}}
\put(100,17){\vector(0,-1){9}}
\put(60,22){\makebox(0,0)[c]{\footnotesize$\gamma_*\alpha$}}
\put(60,3){\makebox(0,0)[c]{\footnotesize$\alpha^\gamma$}}
\put(38,12){\makebox(0,0)[c]{\footnotesize$\gamma$}}
\put(106,12){\makebox(0,0)[c]{\footnotesize$\id\times\gamma$}}
\put(66,14){\makebox(0,0)[c]{\footnotesize$\gamma\circ(\gamma_*\alpha)$}}
\put(69,9){\vector(3,-1){3}}
\dottedline{1}(48,16)(69,9)
\end{picture}

The morphism $\alpha$ is given by a quotient 
$[0\to\mathscr I\to\O_{X_T}\to\O_Z\to 0]$
on $X_T=X_L\times_LT$. Under the identification
$\XHilb{G_L}_LX_L=(\GHilb_KX_K)_L$ the $T$-valued 
point $\alpha$ corresponds to a morphism $T\to\GHilb_KX$ 
of $K$-schemes, that is a quotient 
$[0\to\mathscr I\to\O_{X\times_KT}\to\O_Z\to 0]$ 
on $X\times_KT$, and the structure morphism $f\colon T\to\Spec L$.
We have the correspondences
\[
\begin{array}{rcl}
\alpha&\quad\longleftrightarrow\quad&
\left\{
\begin{array}{l}
[0\to\mathscr I\to\O_{X\times_KT}\to\O_Z\to 0]\\
f\colon T\to\Spec L\\
\end{array}\right.\\
\gamma\circ(\gamma_*\alpha)&\quad\longleftrightarrow
\quad&
\left\{
\begin{array}{l}
[0\to\mathscr I\to\O_{X\times_K\gamma_*T}\to\O_Z\to 0]\\
\gamma\circ(\gamma_*f)\colon\gamma_*T\to\Spec L\\
\end{array}\right.\\
\alpha^\gamma=\gamma\circ(\gamma_*\alpha)\circ\gamma^{-1}&
\quad\longleftrightarrow\quad&
\left\{
\begin{array}{l}
[0\to{\gamma^{-1}}^*\mathscr I\to\O_{X\times_KT}\to
{\gamma^{-1}}^*\O_Z\to 0]\\
f=\gamma\circ(\gamma_*f)\circ\gamma^{-1}\colon T\to\Spec L\\
\end{array}\right.\\

&\quad\longleftrightarrow\quad&
\left\{
\begin{array}{l}
[0\to\gamma_*\mathscr I\to\O_{X\times_KT}\to\O_{\gamma Z}\to 0]\\
f=\gamma\circ(\gamma_*f)\circ\gamma^{-1}\colon T\to\Spec L\\
\end{array}\right.\\
\end{array}
\]
Under the identification $(\GHilb_KX)_L=\XHilb{G_L}_LX_L$
the last morphism corresponds to the quotient
$[0\to\gamma_*\mathscr I\to\O_{X_T}\to\O_{\gamma Z}\to 0]$
on $X_T=X_L\times_LT$. 
\end{proof}

In particular, in the case $X=\A_K^2$ the $\gamma$-conjugate of an 
$L$-valued point given by an ideal $I\sub L[x_1,x_2]$ or a $G_L$-cluster 
$Z\subset\A^2_L$ is given by the $\gamma$-conjugate ideal
$\gamma^{-1}I\subset L[x_1,x_2]$ or the $\gamma$-conjugate 
$G_L$-cluster $\gamma Z\subset\A^2_L$.\\
Every point $x$ of the $L$-scheme $\XHilb{G_L}_L\A^2_L$ such that 
$\kappa(x)=L$ corresponds to a unique $L$-valued point 
$\alpha\colon\Spec L\to\XHilb{G_L}_L\A^2_L$.
The $\gamma$-conjugate point $\gamma x$ corresponds to
the $\gamma$-conjugate $L$-valued point 
$\alpha^\gamma\colon\Spec L\to\XHilb{G_L}_L\A^2_L$.

\begin{cor}\label{cor:conjpoint}
Let $x$ be a closed point of $\XHilb{G_L}_L\A^2_L$ such that
$\kappa(x)=L$, $\alpha\colon\Spec L\to \XHilb{G_L}_L\A^2_L$ the 
corresponding $L$-valued point given by an ideal 
$I\subset L[x_1,x_2]$. 
Then for $\gamma\in\Gamma$ the conjugate point $\gamma x$
corresponds to the $\gamma$-conjugate $L$-valued point
$\alpha^\gamma\colon\Spec L\to\XHilb{G_L}_L\A^2_L$,
which is given by the ideal $\gamma^{-1}I\subset L[x_1,x_2]$.\qed
\end{cor}

\subsection{Conjugate G-sheaves}

Let $G=\Spec A$ be a group scheme over a field $K$, 
$X$ be a $G$-scheme over $K$, let $K\to L$ be a field 
extension and $\Gamma=\Aut_K(L)$. Again, $\Gamma$ operates
on $X_L$ by automorphisms $\gamma\colon X_L\to X_L$ over $K$, 
these are equivariant with respect to the $G$-scheme 
structure of $X_L$ defined in remark \ref{rem:Gsheaves}.(8).

\begin{propdefi} 
Let $\mathscr F$ be a $G_L$-sheaf on $X_L$. 
For $\gamma\in\Gamma$ the $\O_{X_L}$-module $\gamma_*\mathscr F$ has 
a natural $G_L$-sheaf structure given by

\smallskip

\begin{picture}(160,25)
\put(15,0){
\put(40,20){\makebox(0,0)[c]{$\gamma_*s_{X_L}^*\mathscr F$}}
\put(90,20){\makebox(0,0)[c]{$\gamma_*p_{X_L}^*\mathscr F$}}
\put(40,5){\makebox(0,0)[c]{$s_{X_L}^*\gamma_*\mathscr F$}}
\put(90,5){\makebox(0,0)[c]{$p_{X_L}^*\gamma_*\mathscr F$}}
\put(50,20){\vector(1,0){30}}
\put(50,5){\vector(1,0){30}}
\put(40,8){\vector(0,1){9}}
\put(90,8){\vector(0,1){9}}
\put(65,22){\makebox(0,0)[c]{\footnotesize$\gamma_*\lambda^{\mathscr F}$}}
\put(65,7){\makebox(0,0)[c]{\footnotesize$\lambda^{\gamma_*\mathscr F}$}}
\put(39,12){\rotatebox{-90}{\makebox(0,0)[c]{\footnotesize$\sim$}}}
\put(91,12){\rotatebox{90}{\makebox(0,0)[c]{\footnotesize$\sim$}}}
}
\end{picture}

\vspace{-21mm}
\begin{equation}\label{eq:conjGsheaf}\end{equation}

\vspace{9mm}

This $G_L$-sheaf $\gamma_*\mathscr F$ is called the
$\gamma$-conjugate $G_L$-sheaf of $\mathscr F$. 
For a morphism of $G_L$-sheaves $\phi\colon\mathscr F\to\mathscr G$
the morphism $\gamma_*\phi\colon\gamma_*\mathscr F\to\gamma_*\mathscr G$
is a morphism of $G_L$-sheaves between the sheaves 
$\gamma_*\mathscr F$ and $\gamma_*\mathscr G$ with $\gamma$-conjugate 
$G_L$-sheaf structure.\qed
\end{propdefi}

\begin{rem}
This way functors $\gamma_*$ are defined, similarly one may define 
functors $\gamma^*$, then $\gamma_*$ and $(\gamma^{-1})^*$ are 
isomorphic. In the case of trivial operation they preserve
trivial $G$-sheaf structures.
\end{rem}

The functors $\gamma_*$ commute with functors $f_L^*,f_{L*}$ for 
equivariant morphisms $f$ and with bifunctors like $\SHom$ 
and $\otimes$:

\begin{lemma}\label{le:gamma_*}
There are the following natural isomorphisms of $G_L$-sheaves: 

\vspace{-1mm}

\begin{enumerate}[\rm (i)]\setlength{\itemsep}{0mm}
\item For $G_L$-sheaves $\mathscr F,\mathscr G$ on $X_L:\;$ 
$\gamma_*(\mathscr F\otimes_{\O_{X_L}}\mathscr G)\cong
\gamma_*\mathscr F\otimes_{\O_{X_L}}\gamma_*\mathscr G$.
\item Let $f\colon Y\to X$ be an equivariant morphism of $G$-schemes 
over $K$ and $\mathscr F$ a $G_L$-sheaf on $X_L$. 
Then $\gamma_*f_L^*\mathscr F\cong f_L^*\gamma_*\mathscr F$.
\item For quasicoherent $G_L$-sheaves $\mathscr F,\mathscr G$ 
on $X_L$ with $\mathscr F$ finitely presented: 
$\gamma_*\SHom_{\O_{X_L}}(\mathscr F,\mathscr G)\cong
\SHom_{\O_{X_L}}(\gamma_*\mathscr F,\gamma_*\mathscr G)$.
If the $G$-operation on $X$ is trivial, it follows that
$\gamma_*(\SHom^{G_L}_{\O_{X_L}}(\mathscr F,\mathscr G))$
$\cong\SHom^{G_L}_{\O_{X_L}}(\gamma_*\mathscr F,\gamma_*
\mathscr G)$.\qed
\end{enumerate}
\end{lemma}

\begin{rem}\label{rem:F->gamma_*F}
If $\mathscr F\cong\mathscr F'_L$ for some $G$-sheaf 
$\mathscr F'$ on $X$, then there are maps (not $L$-linear) 
$\gamma\colon\mathscr F\to\mathscr F$ resp.\ isomorphisms of 
$G_L$-sheaves $\gamma\colon\mathscr F\to\gamma_*\mathscr F$ on 
$X_L$. For a subsheaf $\mathscr G\sub\mathscr F$
the above isomorphisms of $G_L$-sheaves restrict to 
isomorphisms of $G_L$-sheaves
$\gamma\colon\gamma^{-1}\mathscr G\to\gamma_*\mathscr G$.
\end{rem}

\subsection{Conjugate comodules and representations}

Let $G=\Spec A$ be an affine group scheme over a field $K$,
$X$ be a $G$-scheme over $K$, let $K\to L$ be a field extension
and $\Gamma=\Aut_K(L)$.

\begin{rem}\label{rem:conjHopfalg}
For $\gamma\in\Gamma$ there are maps $\gamma\colon A_L\to A_L$.
Taking the canonically defined conjugate Hopf algebra 
structure on the target, these maps become isomorphisms 
$\gamma\colon A_L\to\gamma_*A_L$ of Hopf algebras over $L$.
They correspond to isomorphisms $\gamma\colon\gamma_*G_L\to G_L$ 
of group schemes over $L$.
\end{rem}

\begin{prop}\label{prop:conjcomod}
Let $\mathscr F$ be a $G_L$-sheaf on $X_L$, $X$ with trivial 
$G$-operation, the $G_L$-sheaf structure equivalent to an 
$A_L$-comodule structure $\rho^{\mathscr F}\colon\mathscr F\to 
A_L\otimes_L\mathscr F$.
Then for $\gamma\in\Gamma$ the $G_L$-sheaf structure of the 
$\gamma$-conjugate $G_L$-sheaf $\gamma_*\mathscr F$ is 
equivalent to the comodule structure $\rho^{\gamma_*\mathscr F}\colon
\gamma_*\mathscr F\to A_L\otimes_L\gamma_*\mathscr F$
determined by commutativity of the diagram

\smallskip

\begin{picture}(160,25)
\put(10,0){
\put(40,20){\makebox(0,0)[c]{$\gamma_*\mathscr F$}}
\put(95,20){\makebox(0,0)[c]{$\gamma_* A_L\otimes_L\gamma_*\mathscr F$}}
\put(40,5){\makebox(0,0)[c]{$\gamma_*\mathscr F$}}
\put(97,5){\makebox(0,0)[c]{$A_L\otimes_L\gamma_*\mathscr F$}}
\put(50,20){\vector(1,0){30}}
\put(50,5){\vector(1,0){30}}
\put(40,8){\vector(0,1){9}}
\put(94.5,8){\vector(0,1){9}}
\put(65,22.5){\makebox(0,0)[c]{\footnotesize$\gamma_*\rho^{\mathscr F}$}}
\put(65,7.5){\makebox(0,0)[c]{\footnotesize$\rho^{\gamma_*\mathscr F}$}}
\put(37,12){\makebox(0,0)[c]{\footnotesize$\id$}}
\put(100,12){\makebox(0,0)[c]{\footnotesize$\gamma\otimes\id$}}
}
\end{picture}

\vspace{-4mm}
\end{prop}

\begin{proof}
Apply the construction mentioned in remark \ref{rem:Gsheaves}.(2) 
to diagram (\ref{eq:conjGsheaf}).
\end{proof}

In the special case of representations the definition
of conjugate $G$-sheaves leads to the notion of a 
conjugate representation:
Instead of a sheaf $\gamma_*\mathscr F$ one has an $L$-vector space
$\gamma_*V$, the vector space structure given by
$(l,v)\mapsto\gamma(l)v$ using the original structure.
The choice of a $K$-structure $V=V'_L$ gives an isomorphism 
$\gamma\colon V\to\gamma_*V$ of $L$-vector spaces and leads to the diagram

\smallskip

\begin{picture}(160,40)
\put(10,0){
\put(40,35){\makebox(0,0)[c]{$\gamma_*V$}}
\put(95,35){\makebox(0,0)[c]{$\gamma_*A_L\otimes_L\gamma_*V$}}
\put(40,20){\makebox(0,0)[c]{$\gamma_*V$}}
\put(97,20){\makebox(0,0)[c]{$A_L\otimes_L\gamma_*V$}}
\put(40,5){\makebox(0,0)[c]{$V$}}
\put(97,5){\makebox(0,0)[c]{$A_L\otimes_LV$}}
\put(50,35){\vector(1,0){30}}
\put(50,20){\vector(1,0){30}}
\put(50,5){\vector(1,0){30}}
\put(40,23){\vector(0,1){9}}
\put(40,8){\vector(0,1){9}}
\put(94.5,23){\vector(0,1){9}}
\put(94.5,8){\vector(0,1){9}}
\put(65,37.5){\makebox(0,0)[c]{\footnotesize$\gamma_*\rho^V$}}
\put(65,22.5){\makebox(0,0)[c]{\footnotesize$\rho^{\gamma_*V}$}}
\put(65,7.5){\makebox(0,0)[c]{\footnotesize$(\rho^V)^\gamma$}}
\put(37,27){\makebox(0,0)[c]{\footnotesize$\id$}}
\put(38,12){\makebox(0,0)[c]{\footnotesize$\gamma$}}
\put(100,27){\makebox(0,0)[c]{\footnotesize$\gamma\otimes\id$}}
\put(100,12){\makebox(0,0)[c]{\footnotesize$\id\otimes\gamma$}}
}
\end{picture}

for definition of the $\gamma$-conjugate $A_L$-comodule structure
$(\rho^V)^\gamma$ on $V$ --- this definition is made, such that
$\gamma\colon (V,(\rho^V)^\gamma)\to(\gamma_*V,\rho^{\gamma_*V})$ is 
an isomorphism of $A_L$-comodules. We write $V^\gamma$
for $V$ with the conjugate $A_L$-comodule structure.

\enlargethispage{5mm}

\begin{rem}\label{rem:gammaU}
Let $V'$ be an $A$-comodule over $K$ and $V=V'_L$. Then as a 
special case of remark \ref{rem:F->gamma_*F} there are
maps $\gamma\colon V\to V$ resp.\ isomorphisms of $A_L$-comodules 
$\gamma\colon V\to\gamma_*V$. For any $A_L$-subcomodule $U\sub V$ 
these restrict to isomorphisms of $A_L$-comodules 
$\gamma^{-1}U\stackrel{\!_\sim}{\longrightarrow}\gamma_*U\cong U^{\gamma}$. 
\end{rem}

\subsection{Decomposition into isotypic components and Galois 
extensions}\label{subsec:decomp-Galoisext}

Let $G=\Spec A$ be an affine group scheme over a field $K$, 
let $K\to L$ be a Galois extension, $\Gamma=\Aut_K(L)$. 
Assume that $A, A_L$ are cosemisimple.\medskip

Recall the relations between the Galois operation on $A_L$
given by maps $\gamma\colon A_L\to A_L$ resp.\ isomorphisms 
$\gamma\colon A_L\to\gamma_*A_L$ of Hopf algebras or of $A_L$-comodules
(see remark \ref{rem:conjHopfalg} or \ref{rem:gammaU}) 
and the decompositions $A=\bigoplus_{j\in J}A_j$ and 
$A_L=\bigoplus_{i\in I}A_{L,i}$ into simple subcoalgebras 
described in corollary \ref{cor:decompcoalg}.
We relate this to conjugation of representations.
The subcoalgebras $A_{L,i}$ are the isotypic components of 
$A_L$ as a left-(or right-)comodule, let $V_i$ be the 
isomorphism class of simple $A_L$-comodules corresponding 
to $A_{L,i}$.
Define an operation of $\Gamma$ on the index set $I$ by 
$V_{\gamma(i)}=V_i^\gamma$. Using remark \ref{rem:gammaU}
one obtains:

\begin{lemma}\label{le:gammaA_L,i}
$\gamma^{-1}A_{L,i}=A_{L,\gamma(i)}$.\qed
\end{lemma}

The decomposition of $A$ into simple subcoalgebras 
$A=\bigoplus_jA_j$ gives decompositions of representations 
and more generally of $G$-sheaves on $G$-schemes with 
trivial $G$-operation into isotypic components corresponding 
to the $A_j$ (see remark \ref{rem:Gsheaves}.(3)).
After base extension one has decompositions of 
$G_L$-sheaves, we compare it with the decompositions 
coming from the decomposition of $A_L$ into simple 
subcoalgebras.

\begin{prop}\label{prop:conjdecomp}
Let $X$ be a $G$-scheme with trivial operation, $\mathscr F$ 
a $G$-sheaf on $X$ and let 
\[\textstyle\mathscr F=\bigoplus_j\mathscr F_j,\qquad 
\mathscr F_L=\bigoplus_i\mathscr F_{L,i}\]
be the decompositions into isotypic components as a $G$-sheaf 
resp.\ $G_L$-sheaf. Then:

\vspace{-2mm}

\begin{enumerate}[\rm (i)]\setlength{\itemsep}{0.5mm}
\item $\mathscr F_L=\bigoplus_i\mathscr F_{L,i}$ is a
refinement of $\mathscr F_L=\bigoplus_j(\mathscr F_j)_L$.
\item The operation of $\Gamma$ on $\mathscr F_L$ (see remark 
\ref{rem:F->gamma_*F}) permutes the isotypic components
$\mathscr F_{L,i}$ of $\mathscr F_L$. 
It is $\gamma^{-1}\mathscr F_{L,i}=\mathscr F_{L,\gamma(i)}$,
if $V_{\gamma(i)}=V_i^\gamma$.
\item $(\mathscr F_j)_L=\sum_{\gamma\in\Gamma}
\gamma\mathscr F_{L,i}$, if $\mathscr F_{L,i}$ is a summand
of $(\mathscr F_j)_L$.
\end{enumerate}
\end{prop}
\begin{proof}[Sketch of proof.] Combine remark \ref{rem:F->gamma_*F}, 
proposition \ref{prop:conjcomod} and lemma \ref{le:gammaA_L,i} 
with corollary \ref{cor:decompcoalg}.
\end{proof}

\begin{cor}\label{cor:irredrep}
$\Gamma$ operates by $V_i\mapsto V_i^\gamma$ 
on the set $\{V_i\:|\:i\in I\}$ of isomorphism classes 
of irreducible representations of $G_L$. 
The subsets of $\{V_i\:|\:i\in I\}$, which occur by decomposing 
irreducible representations of $G$ over $K$ as representations 
over $L$, are exactly the $\Gamma$-orbits.\qed
\end{cor}

For similar results in the representation theory of finite groups
see e.g.~\cite[Vol. I, \S 7B]{CR}.

\section{The finite subgroup schemes of $\SL(2,K)$: representations
and graphs}

In this section $K$ denotes a field of characteristic $0$.

\subsection{The finite subgroups of $\SL(2,C)$}

By the well known classification any finite subgroup 
$G\subset\SL(2,C)$, $C$ an algebraically closed field 
of characteristic $0$, is isomorphic to one of the 
following groups (presentations and character tables are 
listed in appendix \ref{app:charsubgrSL2}).\smallskip

$\Z/n\Z$ (cyclic group of order $n$), $n\geq 1$\\ 
$\BD_n$ (binary dihedral group of order $4n$), $n\geq 2$\\ 
$\BT$ (binary tetrahedral group)\\
$\BO$ (binary octahedral group)\\
$\BI$ (binary icosahedral group).

\subsection{Representation graphs}

In the following definition we will introduce the (extended)
representation graph as an in general directed graph.
A loop is defined to be an edge emanating from and terminating at 
the same vertex. 
In addition we will attach a natural number called multiplicity
to any vertex, and for homomorphisms of graphs in addition 
we will require, that for any vertex of the target its multiplicity 
is the sum of the multiplicities of its preimages.

\begin{defi}
The extended representation graph $\graph(G,V)$ associated to 
a finite subgroup scheme $G$ of $\GL(n,K)$, $V$ the given 
$n$-dimensional representation, is defined as the following
directed graph:\smallskip

- vertices. A vertex of multiplicity $n$ for each irreducible
representation of $G$ over $K$ which decomposes over the algebraic 
closure of $K$ into $n$ irreducible representations.\smallskip

- edges. Vertices $V_i$ and $V_j$ are connected by
$\dim_K\Hom^G_K(V_i,V\otimes_K V_j)$ directed edges from $V_i$
to $V_j$. In particular any vertex $V_i$ has
$\dim_K\Hom^G_K(V_i,V\otimes_K V_i)$ directed loops.\smallskip

Define the representation graph to be the graph, which arises 
by leaving out the trivial representation and all edges emanating
from or terminating at the trivial representation.
\end{defi}

We say that a graph is undirected, if between any two different
vertices the numbers of directed edges of both directions coincide
and for any vertex the number of directed loops is even.

Then one can form a graph having only undirected edges by
defining {\it(number of undirected edges between $V_i$ and $V_j$) 
$:=$ (number of directed edges from $V_i$ to $V_j$) $=$ (number of 
directed edges from $V_j$ to $V_i$)} for different vertices $V_i,V_j$ 
and {\it(number of undirected loops of $V_i$) $:=$ $\frac{1}{2}$(number 
of directed loops of $V_i$)} for any vertex $V_i$.

\begin{rem}\label{rem:reprgraph}\ \\
(1) For $G\subset\SL(2,K)$ the (extended) representation graph 
is undirected. There is the isomorphism
$\Hom^G_K(V_i\otimes_K V,V_j)\cong\Hom^G_K(V_i,V\otimes_K V_j)$,
which follows from the isomorphism 
$\Hom^G_K(V_i\otimes_K V,V_j)\cong\Hom^G_K(V_i,V^\vee\otimes_K V_j)$ 
and the fact that the $2$-dimensional representation $V$ given by 
inclusion $G\to\SL(2,K)$ is self-dual.
Further, that the number of directed loops of any vertex is even, 
follows from the fact that over the algebraic closure $C$ one has 
$\dim_C\Hom^G_C(U_i,V_C\otimes_C U_i)=0$ for irreducible $U_i$ 
over $C$.\\
(2) There is a definition of (extended) representation graph with 
another description of the edges: vertices $V_i$ and $V_j$ are 
connected by $a_{ij}$ edges from $V_i$ to $V_j$, where 
$V\otimes_K V_j=a_{ij}V_i\oplus\textit{other summands}$.
The two definitions coincide over algebraically closed fields,
always one has $a_{ij}\leq\dim_K\Hom^G_K(V_i,V\otimes_KV_j)$,
inequality comes from the presence of nontrivial automorphisms.
\end{rem}

\begin{defi}
For a finite subgroup scheme $G\subset\SL(2,K)$, $V$ the given
$2$-dimensional representation, define a $\Z$-bilinear form
$\<\cdot,\cdot\>$ on the representation ring of $G$ by
\[\<V_i,V_j\>:=\dim_K\Hom^G_K(V_i,V\otimes_KV_j)-
2\dim_K\Hom^G_K(V_i,V_j)\]
\end{defi}

\begin{rem} The form $\<\cdot,\cdot\>$ determines and is 
determined by the extended representation graph
(the second equation follows from the fact, that 
$\dim_K\Hom^G_K(V_i,V_i)=\textit{multiplicity of $V_i$}$): 
\[\begin{array}{l}
\<V_i,V_j\>=\<V_j,V_i\>=\textit{number of undirected edges between 
$V_i$ and $V_j$},\;\;\textit{if}\; V_i\not\cong V_j\\
\frac{1}{2}\<V_i,V_i\>=\textit{number of undirected loops of $V_i$}
-\textit{multiplicity of $V_i$}\\
\end{array}
\]
\end{rem}

\subsection{Representation graphs and field extensions}

Let $K\to L$ be a Galois extension, $\Gamma=\Aut_K(L)$ and let
$G$ be a finite subgroup scheme of $\SL(2,K)$.\medskip

An irreducible representation $W$ of $G$ over $K$ 
decomposes as a representation of $G_L$ over $L$
into isotypic components $W=\bigoplus_i U_i$ which are 
$\Gamma$-conjugate by proposition \ref{prop:conjdecomp}.
Every $U_i$ decomposes into irreducible components
$U_i=V_i^{\oplus m}$ (the same $m$ for all $i$ because of 
$\Gamma$-conjugacy). In the following we will write 
$m(W,L/K)$ for this number. It is related to the Schur index
in the representation theory of finite groups 
(see e.g.~\cite[Vol. II, \S 74]{CR}).

\begin{prop}\label{prop:m}
For finite subgroup schemes $G$ of $\SL(2,K)$ it is
$m(W_j,L/K)=1$ for every irreducible representation 
$W_j$ of $G$. It follows that $W_j$ decomposes over $L$ 
into a direct sum $(W_j)_L\cong\bigoplus_iV_i$ of 
$\gamma$-conjugate irreducible representations $V_i$ of $G_L$ 
nonisomorphic to each other.
\end{prop}
\begin{proof} We may assume $L$ algebraically closed.
Further we may assume that $G$ is not cyclic.

The natural $2$-dimensional representation $W$ given by inclusion
$G\subset\SL(2,K)$ does satisfy $m(W,L/K)=1$, because it is
irreducible over $L$.

Following the discussion below without using this proposition
one obtains the graphs in section 4.4 without multiplicities
of vertices and edges but one knows which vertices over the algebraic 
closure may form a vertex over $K$ and which vertices are connected.
Argue that if an irreducible representation $W_i$ satisfies 
$m(W_i,L/K)=1$, then any irreducible $W_j$ connected to $W_i$ 
in the representation graph has to satisfy this property as well.
\end{proof}

There is a morphism of graphs $\graph(G_L,W_L)\to\graph(G,W)$
(resp.\ of the nonextended graphs, the following applies 
to them as well).
For $W_j$ an irreducible representation of $G$ the base extension
$(W_j)_L$ is a sum $(W_j)_L=\bigoplus_{i}V_i$ of irreducible 
representations of $G_L$ nonisomorphic to each other by proposition 
\ref{prop:m}. The morphism $\graph(G_L,W_L)\to\graph(G,W)$ maps 
components of $(W_j)_L$ to $W_j$, thereby their multiplicities are 
added. 
Further, for irreducible representations $W_j,W_{j'}$ of $G$ there is 
a bijection between the set of edges between $W_j$ and $W_{j'}$ and 
the union of the sets of edges between the irreducible components of 
$(W_j)_L$ and $(W_{j'})_L$, again using proposition \ref{prop:m}
$(W_j)_L$ and $(W_{j'})_L$ are sums $(W_j)_L=\bigoplus_{i}V_i$, 
$(W_{j'})_L=\bigoplus_{i'}V_{i'}$ of irreducible representations 
of $G_L$ nonisomorphic to each other and one has 
\[
\begin{array}{rcl}
\dim_K\Hom^G_K(W_j\otimes_K W,W_{j'})&=&
\dim_L(\Hom^G_K(W_j\otimes_K W,W_{j'})\otimes_KL)\\
&=&\dim_L\Hom^G_L((W_j)_L\otimes_L W_L,(W_{j'})_L)\\
&=&\dim_L\Hom^G_L(\bigoplus_iV_i\otimes_L W_L,\bigoplus_{i'}V_{i'})\\
&=&\sum_{i,i'}\dim_L\Hom^G_L(V_i\otimes_L W_L,V_{i'})\\
\end{array}
\]
The Galois group $\Gamma$ operates on $\graph(G_L,W_L)$ by graph 
automorphisms. Irreducible representations are mapped to conjugate 
representations and equivariant homomorphisms to the conjugate 
homomorphisms. The vertices of $\graph(G,W)$ correspond to 
$\Gamma$-orbits of vertices of $\graph(G_L,W_L)$ by corollary 
\ref{cor:irredrep}.

\begin{prop}\label{prop:repgraph}
The (extended) representation graph of $G$ arises by identifying 
the elements of $\Gamma$-orbits of vertices of the (extended) 
representation graph of $G_L$, adding multiplicities. 
The edges between vertices $W_j$ and $W_{j'}$ are in bijection 
with the edges between the isomorphism classes of irreducible 
components of $(W_j)_L$ and $(W_{j'})_L$.\qed
\end{prop}

\subsection{The representation graphs of the finite subgroup schemes
of $\SL(2,K)$}\label{subsec:repgraphs}

As extended representation graph of a finite subgroup 
scheme of $\SL(2,K)$ with respect to the natural $2$-dimensional 
representation the following graphs can occur.
We list the extended representation graphs $\graph(G,V)$
of the finite subgroups of $\SL(2,C)$ for $C$ algebraically closed, 
their groups of automorphisms leaving the trivial representation 
fixed and the possible extended representation graphs for finite 
subgroup schemes over non algebraically closed $K$, which after 
suitable base extension become the graph $\graph(G,V)$.
We use the symbol $\circ$ for the trivial representation.\bigskip

\pagebreak
- Cyclic groups

\begin{picture}(160,67)

\put(0,35){
\put(7,25){\makebox(0,0)[l]{$(A_{2n}), n\geq 1$}}
\put(12,15){\makebox(0,0)[l]{$\Z/2\Z$}}
\put(32,15){
\filltype{white}
\put(0,0){\circle*{1.2}}
\filltype{black}
\put(10,6){\circle*{1.2}}
\put(10,-6){\circle*{1.2}}
\put(36,6){\circle*{1.2}}
\put(36,-6){\circle*{1.2}}
\put(46,6){\circle*{1.2}}
\put(46,-6){\circle*{1.2}}
\drawline(1,0.6)(9,5.4)
\drawline(1,-0.6)(9,-5.4)
\drawline(11,6)(15,6)
\drawline(11,-6)(15,-6)
\dottedline{1}(15,6)(31,6)
\dottedline{1}(15,-6)(31,-6)
\drawline(31,6)(35,6)
\drawline(31,-6)(35,-6)
\drawline(37,6)(45,6)
\drawline(37,-6)(45,-6)
\drawline(46,5)(46,-5)
}
\put(90,25){\makebox(0,0)[l]{$(A_{2n})'$}}
\put(100,15){
\filltype{white}
\put(0,0){\circle*{1.2}}
\put(47.5,0){\circle*{3}}
\filltype{black}
\put(10,0.6){\circle*{1.2}}
\put(10,-0.6){\circle*{1.2}}
\put(36,0.6){\circle*{1.2}}
\put(36,-0.6){\circle*{1.2}}
\put(46,0.6){\circle*{1.2}}
\put(46,-0.6){\circle*{1.2}}
\drawline(1,0.6)(9,0.6)
\drawline(1,-0.6)(9,-0.6)
\drawline(11,0.6)(15,0.6)
\drawline(11,-0.6)(15,-0.6)
\dottedline{1}(15,0.6)(31,0.6)
\dottedline{1}(15,-0.6)(31,-0.6)
\drawline(31,0.6)(35,0.6)
\drawline(31,-0.6)(35,-0.6)
\drawline(37,0.6)(45,0.6)
\drawline(37,-0.6)(45,-0.6)
}
}

\put(0,10){
\put(7,25){\makebox(0,0)[l]{$(A_{2n+1}), n\geq 1$}}
\put(12,15){\makebox(0,0)[l]{$\Z/2\Z$}}
\put(32,15){
\filltype{white}
\put(0,0){\circle*{1.2}}
\filltype{black}
\put(10,6){\circle*{1.2}}
\put(10,-6){\circle*{1.2}}
\put(36,6){\circle*{1.2}}
\put(36,-6){\circle*{1.2}}
\put(46,0){\circle*{1.2}}
\drawline(1,0.6)(9,5.4)
\drawline(1,-0.6)(9,-5.4)
\drawline(11,6)(15,6)
\drawline(11,-6)(15,-6)
\dottedline{1}(15,6)(31,6)
\dottedline{1}(15,-6)(31,-6)
\drawline(31,6)(35,6)
\drawline(31,-6)(35,-6)
\drawline(37,5.4)(45,0.6)
\drawline(37,-5.4)(45,-0.6)
}
\put(90,25){\makebox(0,0)[l]{$(A_{2n+1})'$}}
\put(100,15){
\filltype{white}
\put(0,0){\circle*{1.2}}
\filltype{black}
\put(10,0.6){\circle*{1.2}}
\put(10,-0.6){\circle*{1.2}}
\put(36,0.6){\circle*{1.2}}
\put(36,-0.6){\circle*{1.2}}
\put(46,0){\circle*{1.2}}
\drawline(1,0.6)(9,0.6)
\drawline(1,-0.6)(9,-0.6)
\drawline(11,0.6)(15,0.6)
\drawline(11,-0.6)(15,-0.6)
\dottedline{1}(15,0.6)(31,0.6)
\dottedline{1}(15,-0.6)(31,-0.6)
\drawline(31,0.6)(35,0.6)
\drawline(31,-0.6)(35,-0.6)
\drawline(37,0.6)(45,0.6)
\drawline(37,-0.6)(45,-0.6)
}
}

\put(7,13){\makebox(0,0)[l]{$(A_1)$}}
\put(12,5){\makebox(0,0)[l]{$\{\id\}$}}
\put(32,5){
\filltype{white}
\put(0,0){\circle*{1.2}}
\filltype{black}
\put(10,0){\circle*{1.2}}
\drawline(1,0.6)(9,0.6)
\drawline(1,-0.6)(9,-0.6)
}
\end{picture}

\medskip

- Binary dihedral groups

\begin{picture}(160,52)

\put(0,20){
\put(7,25){\makebox(0,0)[l]{$(D_n), n\geq 5$}}
\put(12,15){\makebox(0,0)[l]{$\Z/2\Z$}}
\put(32,15){
\filltype{white}
\put(0,6){\circle*{1.2}}
\filltype{black}
\put(0,-6){\circle*{1.2}}
\put(10,0){\circle*{1.2}}
\put(36,0){\circle*{1.2}}
\put(46,6){\circle*{1.2}}
\put(46,-6){\circle*{1.2}}
\drawline(1,5.4)(9,0.6)
\drawline(1,-5.4)(9,-0.6)
\drawline(11,0)(15,0)
\dottedline{1}(15,0)(31,0)
\drawline(31,0)(35,0)
\drawline(37,0.6)(45,5.4)
\drawline(37,-0.6)(45,-5.4)
}
\put(90,25){\makebox(0,0)[l]{$(D_n)'$}}
\put(100,15){
\filltype{white}
\put(0,6){\circle*{1.2}}
\filltype{black}
\put(0,-6){\circle*{1.2}}
\put(10,0){\circle*{1.2}}
\put(36,0){\circle*{1.2}}
\put(46,0.6){\circle*{1.2}}
\put(46,-0.6){\circle*{1.2}}
\drawline(1,5.4)(9,0.6)
\drawline(1,-5.4)(9,-0.6)
\drawline(11,0)(15,0)
\dottedline{1}(15,0)(31,0)
\drawline(31,0)(35,0)
\drawline(37,0.6)(45,0.6)
\drawline(37,-0.6)(45,-0.6)
}
}

\put(0,-5){
\put(7,25){\makebox(0,0)[l]{$(D_4)$}}
\put(12,15){\makebox(0,0)[l]{$S_3$}}
\put(32,15){
\filltype{white}
\put(0,6){\circle*{1.2}}
\filltype{black}
\put(0,-6){\circle*{1.2}}
\put(10,0){\circle*{1.2}}
\put(20,6){\circle*{1.2}}
\put(20,-6){\circle*{1.2}}
\drawline(1,5.4)(9,0.6)
\drawline(1,-5.4)(9,-0.6)
\drawline(11,0.6)(19,5.4)
\drawline(11,-0.6)(19,-5.4)
}
\put(75,25){\makebox(0,0)[l]{$(D_4)'$}}
\put(90,15){
\filltype{white}
\put(0,6){\circle*{1.2}}
\filltype{black}
\put(0,-6){\circle*{1.2}}
\put(10,0){\circle*{1.2}}
\put(20,0.6){\circle*{1.2}}
\put(20,-0.6){\circle*{1.2}}
\drawline(1,5.4)(9,0.6)
\drawline(1,-5.4)(9,-0.6)
\drawline(11,0.6)(19,0.6)
\drawline(11,-0.6)(19,-0.6)
}
\put(120,25){\makebox(0,0)[l]{$(D_4)''$}}
\put(126,15){
\filltype{white}
\put(0,0){\circle*{1.2}}
\filltype{black}
\put(10,0){\circle*{1.2}}
\put(20,1){\circle*{1.2}}
\put(20,0){\circle*{1.2}}
\put(20,-1){\circle*{1.2}}
\drawline(1,0)(9,0)
\drawline(11,0.8)(19,0.8)
\drawline(11,0)(19,0)
\drawline(11,-0.8)(19,-0.8)
}
}
\end{picture}

\smallskip

- Binary tetrahedral group

\begin{picture}(160,32)
\put(0,0){
\put(7,25){\makebox(0,0)[l]{$(E_6)$}}
\put(12,15){\makebox(0,0)[l]{$\Z/2\Z$}}
\put(32,15){
\filltype{white}
\put(0,0){\circle*{1.2}}
\filltype{black}
\put(10,0){\circle*{1.2}}
\put(20,0){\circle*{1.2}}
\put(30,6){\circle*{1.2}}
\put(30,-6){\circle*{1.2}}
\put(40,6){\circle*{1.2}}
\put(40,-6){\circle*{1.2}}
\drawline(1,0)(9,0)
\drawline(11,0)(19,0)
\drawline(21,0.6)(29,5.4)
\drawline(21,-0.6)(29,-5.4)
\drawline(31,6)(39,6)
\drawline(31,-6)(39,-6)
}
\put(90,25){\makebox(0,0)[l]{$(E_6)'$}}
\put(100,15){
\filltype{white}
\put(0,0){\circle*{1.2}}
\filltype{black}
\put(10,0){\circle*{1.2}}
\put(20,0){\circle*{1.2}}
\put(30,0.6){\circle*{1.2}}
\put(30,-0.6){\circle*{1.2}}
\put(40,0.6){\circle*{1.2}}
\put(40,-0.6){\circle*{1.2}}
\drawline(1,0)(9,0)
\drawline(11,0)(19,0)
\drawline(21,0.6)(29,0.6)
\drawline(21,-0.6)(29,-0.6)
\drawline(31,0.6)(39,0.6)
\drawline(31,-0.6)(39,-0.6)
}
}
\end{picture}

- Binary octahedral group

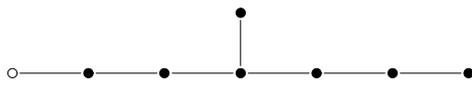
\begin{picture}(160,24)
\put(0,-8){
\put(7,25){\makebox(0,0)[l]{$(E_7)$}}
\put(12,15){\makebox(0,0)[l]{$\{\id\}$}}
\put(32,15){
\filltype{white}
\put(0,0){\circle*{1.2}}
\filltype{black}
\put(10,0){\circle*{1.2}}
\put(20,0){\circle*{1.2}}
\put(30,0){\circle*{1.2}}
\put(40,0){\circle*{1.2}}
\put(50,0){\circle*{1.2}}
\put(60,0){\circle*{1.2}}
\put(30,8){\circle*{1.2}}
\drawline(1,0)(9,0)
\drawline(11,0)(19,0)
\drawline(21,0)(29,0)
\drawline(31,0)(39,0)
\drawline(41,0)(49,0)
\drawline(51,0)(59,0)
\drawline(30,1)(30,7)
}
}
\end{picture}

- Binary icosahedral group

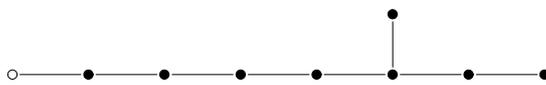
\begin{picture}(160,24)
\put(0,-8){
\put(7,25){\makebox(0,0)[l]{$(E_8)$}}
\put(12,15){\makebox(0,0)[l]{$\{\id\}$}}
\put(32,15){
\filltype{white}
\put(0,0){\circle*{1.2}}
\filltype{black}
\put(10,0){\circle*{1.2}}
\put(20,0){\circle*{1.2}}
\put(30,0){\circle*{1.2}}
\put(40,0){\circle*{1.2}}
\put(50,0){\circle*{1.2}}
\put(60,0){\circle*{1.2}}
\put(70,0){\circle*{1.2}}
\put(50,8){\circle*{1.2}}
\drawline(1,0)(9,0)
\drawline(11,0)(19,0)
\drawline(21,0)(29,0)
\drawline(31,0)(39,0)
\drawline(41,0)(49,0)
\drawline(51,0)(59,0)
\drawline(61,0)(69,0)
\drawline(50,1)(50,7)
}
}
\end{picture}

\begin{rem}
Taking $\frac{2}{-\<V,V\>}V$ for the isomorphism classes 
of irreducible representations $V$ as simple roots one can form 
the Dynkin diagram with respect to the form $-\<\cdot,\cdot\>$ 
(see e.g.~\cite[Groupes et alg\`ebres de Lie]{B}). 
Between (extended) representation graphs and (extended) Dynkin 
diagrams there is the correspondence
\[
\begin{array}{ccccccccccc}
(A_n)&(A_2)'&(A_{2n+1})'&(A_{2n+2})'&(D_n)&(D_n)'&(D_4)''&(E_6)&(E_6)'&
(E_7)&(E_8)\\
(A_n)&(C_1)=(A_1)&(C_{n+1})&(C_{n+1})&(D_n)&(B_{n-1})&(G_2)&(E_6)&(F_4)
&(E_7)&(E_8)\\
\end{array}
\]
A long time ago, the occurrence of the remaining Dynkin diagrams
of types $(B_n)$, $(C_n)$, $(F_4)$, $(G_2)$ as resolution graphs
had been observed in \cite{Li69} with a slightly different 
assignment of the non extended diagrams to the resolutions of 
these singularities, see also \cite{Sl}.
\end{rem}

\subsection{Finite subgroups of $\SL(2,K)$}

Given a field $K$ of characteristic $0$, it is a natural question,
which of the finite subgroups $G\subset\SL(2,C)$, $C$ the 
algebraic closure of $K$, are realisable over the subfield $K$ 
as subgroups (not just as subgroup schemes), that is, there is 
an injective representation of the group $G$ in $\SL(2,K)$.\\
For a finite subgroup $G$ of $\SL(2,C)$ to occur as a subgroup of 
$\SL(2,K)$ it is necessary and sufficient that the given
$2$-dimensional representation in $\SL(2,C)$ is realisable 
over $K$. This is easy to show using the classification and 
the irreducible representations of the individual groups.
If a representation of a group $G$ over $C$ is realisable 
over $K$, necessarily its character has values in $K$. 
For the finite subgroups of $\SL(2,C)$ and the natural 
representation given by inclusion this means:\smallskip

$\Z/n\Z$: $\xi+\xi^{-1}\in K$, \ $\xi\in C$ a primitive $n$-th 
root of unity.\\
$\BD_n$: $\xi+\xi^{-1}\in K$, \ $\xi\in C$ a primitive $2n$-th 
root of unity. \\
$\BT$: no condition.\ \ $\BO$: $\sqrt{2}\in K$.\ \
$\BI$: $\sqrt{5}\in K$.\smallskip

To formulate sufficient conditions, we introduce the following 
notation:

\begin{defi} {\rm (\cite[Part I, Chapter III, \S 1]{Se})}. 
For a field $K$ the Hilbert symbol $\HSymb{\cdot}{\cdot}_K$ 
is the map $K^*\times K^*\to\{-1,1\}$ defined by 
$\HSymb{a}{b}_K=1$, if the equation $z^2-ax^2-by^2=0$ 
has a solution $(x,y,z)\in K^3\setminus\{(0,0,0)\}$, 
and $\HSymb{a}{b}_K=-1$ otherwise.
\end{defi}

\begin{rem}
It is $\HSymb{-1}{b}_K=1$ if and only if $x^2-by^2=-1$ has a 
solution $(x,y)\in K^2$.
\end{rem}

\begin{thm}\label{thm:subgroupsSL(2,K)}
Let $G$ be a finite subgroup of $\SL(2,C)$ such that the values 
of the character of the natural representation given by inclusion
are contained in $K$. Then:

\vspace{-2.5mm}

\begin{enumerate}[\rm (i)]\setlength{\itemsep}{0mm}
\item If $G\cong\Z/n\Z$, then $G$ is realisable over $K$.
\item If $G\cong \BD_n$, let $\xi\in C$ be a primitive $2n$-th 
root of unity and $c:=\frac{1}{2}(\xi+\xi^{-1})$. 
Then $G$ is isomorphic to a subgroup of $\SL(2,K)$ if and only if 
$\HSymb{-1}{c^2-1}_K=1$. 
\item If $G\cong \BT,\BO$ or $\BI$, then $G$ is isomorphic to a 
subgroup of $\SL(2,K)$ if and only if $\HSymb{-1}{-1}_K=1$.
\end{enumerate}
\end{thm}
\begin{proof}
(i) For $n\geq 3$ let $\xi$ be a primitive $n$-th root of unity
and $c:=\frac{1}{2}(\xi+\xi^{-1})$. By assumption $c\in K$.
Then $\Z/n\Z$ is realisable over $K$, there is the representation
\[
\Z/n\Z\to\SL(2,K),\qquad \overline{1}\mapsto
\left(\begin{matrix}0&-1\\1&2c\end{matrix}\right)
\]
(ii) Let $G=BD_n=\<\sigma,\tau\:|\:\tau^2=\sigma^n=(\tau\sigma)^2\>$
(then the element $\tau^2=\sigma^n=(\tau\sigma)^2$ has order $2$) 
and let $\xi$ be a primitive $2n$-th root of unity. 
Then $G$ is realisable as a subgroup of $\SL(2,K)$ if and only if
the representation given by
\begin{equation}\label{eq1}
\sigma\mapsto\left(\begin{matrix}\xi&0\\0&\xi^{-1}
\end{matrix}\right),\qquad
\tau\mapsto\left(\begin{matrix}0&-1\\1&0\end{matrix}\right)
\end{equation}
is realisable over $K$.\smallskip

{\it The representation {\rm (\ref{eq1})} is realisable over $K$ 
if and only if there is a $2\times 2$-matrix $M_\tau$ 
over $K$ having the properties}
\begin{equation}\label{eq2}\textstyle
\det(M_\tau)=1,\;\;\ord(M_\tau)=4,\;\;
(M_\tau M_\sigma)^2=-\mathds 1,\qquad
\textit{where}\; M_\sigma=\left(\begin{smallmatrix}0&-1\\1&2c\end{smallmatrix}
\right),\;\; c=\frac{1}{2}(\xi+\xi^{-1}).
\end{equation}
If the representation (\ref{eq1}) is realisable over $K$,
then with respect to a suitable basis it maps 
$\sigma\mapsto M_\sigma$ and the image of $\tau$ is a matrix
satisfying the properties (\ref{eq2}).\\
On the other hand, if $M_\tau$ is a matrix having these properties,
then $\sigma\mapsto M_\sigma$, $\tau\mapsto M_\tau$ is a
representation of $G$ in $\SL(2,K)$, which is easily seen to be
isomorphic to the representation (\ref{eq1}).\smallskip

{\it There is a $2\times 2$-matrix $M_\tau$ over $K$ having the 
properties {\rm (\ref{eq2})} if and only if the equation
\begin{equation}\label{eq3}
x^2+y^2-2cxy+1=0
\end{equation}
has a solution $(x,y)\in K^2$.}

A matrix  $M_\tau=\left(\begin{smallmatrix}\alpha&\beta\\
\gamma&\delta\end{smallmatrix}\right)$ satisfies the 
conditions (\ref{eq2}) if and only if 
$(\alpha,\beta,\gamma,\delta)\in K^4$ is a solution of
$\alpha\delta-\beta\gamma-1=0,\quad \alpha+\delta=0,\quad
\beta+2c\delta-\gamma=0$. Such an element of $K^4$ exists if
and only if there exists a solution $(x,y)\in K^2$ 
of equation (\ref{eq3}).\smallskip

{\it The equation {\rm (\ref{eq3})} has a solution 
$(x,y)\in K^2$ if and only if $\HSymb{-1}{c^2-1}_K=1$.}\smallskip

We write the equation $x^2+y^2-2cxy+1=0$
as $(x,y)\left(\begin{smallmatrix}1&-c\\-c&1
\end{smallmatrix}\right)\left(\begin{smallmatrix}x\\y
\end{smallmatrix}\right)=-1$. After diagonalisation 
$(x,y)\left(\begin{smallmatrix}1&0\\0&1-c^2
\end{smallmatrix}\right)
\left(\begin{smallmatrix}x\\y\end{smallmatrix}\right)=-1$
or $x^2+(1-c^2)y^2+1=0$. This equation has a solution $(x,y)\in K^2$ 
if and only if $\HSymb{-1}{c^2-1}_K=1$.\smallskip

(iii) Let $G=\BT, \BO$ or $\BI$, that is 
$G=\<a,b\:|\:a^3=b^k=(ab)^2\>$ for $k\in\{3,4,5\}$.
Let $\xi$ be a primitive $2k$-th root of unity and 
$c=\frac{1}{2}(\xi+\xi^{-1})$. As in (ii), using the subgroup 
$\<b\>$ instead of $\<\sigma\>$, we obtain:\smallskip
  
{\it $G$ is isomorphic to a subgroup of $\SL(2,K)$ if and only if
there is a solution $(x,y)\in K^2$ of the equation} 
\begin{equation}\label{eq4}
x^2+y^2-2cxy-x+2cy+1=0
\end{equation}

Next we show:\smallskip

{\it Equation {\rm(\ref{eq4})} has a solution $(x,y)\in K^2$
if and only if $\HSymb{-1}{(2c)^2-3}_K=1$.}

Equation (\ref{eq4}) has a solution if and only if
$(x,y,z)\left(\begin{smallmatrix}
1&-c&-1/2\\-c&1&c\\-1/2&c&1\\
\end{smallmatrix}\right)
\left(\begin{smallmatrix}x\\y\\z\end{smallmatrix}\right)=0$
has a solution $(x,y,z)\in K^3$ with $z\neq 0$.
The existence of a solution with $z\neq 0$ is equivalent to 
the existence of a solution $(x,y,z)\in K^3\setminus\{(0,0,0)\}$ 
(if $(x,y,0)$ is a solution, then $(x,y,x-2cy)$ as well).
After diagonalisation: $(x,y,z)\left(\begin{smallmatrix}
1&0&0\\0&1&0\\0&0&3-(2c)^2\\
\end{smallmatrix}\right)
\left(\begin{smallmatrix}x\\y\\z\end{smallmatrix}\right)=0$.
The existence of a solution $(x,y,z)\in K^3\setminus\{(0,0,0)\}$
for this equation is equivalent to $\HSymb{-1}{(2c)^2-3}_K=1$.\smallskip

For the individual groups we obtain:\\
\begin{tabular}{ll}
$\BT$:&$c=\frac{1}{2}$, $\HSymb{-1}{-2}_K=1$.\\
$\BO$:&$c=\frac{1}{\sqrt{2}}$, $\HSymb{-1}{-1}_K=1$.\\
$\BI$:&$c=\frac{1}{4}(1\pm\sqrt{5})$, 
$\HSymb{-1}{\frac{1}{2}(-3\pm\sqrt{5})}_K=1$.\\
\end{tabular}\\
Each of these conditions is equivalent to $\HSymb{-1}{-1}_K=1$.
For $BI$: $\frac{1}{2}(3\pm\sqrt{5})=(\frac{1}{2}(1\pm\sqrt{5}))^2$.
For $BT$ one has maps between solutions 
$(x,y)$ for $x^2+y^2=-1$ corresponding to $\HSymb{-1}{-1}_K$
and $(x',y')$ for ${x'}^2+2{y'}^2=-1$ corresponding to 
$\HSymb{-1}{-2}_K$ given by 
$x=\frac{x'+1}{2y'}\leftrightarrow x'=\frac{x+y}{x-y}$, 
$y=\frac{x'-1}{2y'}\leftrightarrow y'=\frac{1}{x-y}$ 
for $x\neq y$ resp.\ $y'\neq 0$
and by $(x,x)\mapsto (0,x)$, $(x',0)\mapsfrom(x',0)$.
\end{proof}

\section{McKay correspondence for $G\subset\SL(2,K)$}

Let $G$ be a finite subgroup scheme of $\SL(2,K)$,
$K$ a field of characteristic $0$, and $C$ the algebraic closure 
of $K$.
There is the geometric quotient $\pi\colon\A^2_K\to\A^2_K/G$ 
and the natural morphism $\tau\colon\GHilb_K\A_K^2\to\A_K^2/G$,
which is the minimal resolution of this quotient singularity.

\subsection{The exceptional divisor and the intersection graph}

Define the exceptional divisor $E$ by
\[E:=\tau^{-1}(\overline{O})\]
where $\overline{O}=\pi(O)$, $O$ the origin of $\A^2_K$.
In general $E$ is not reduced, denote by $E_{\rm red}$
the underlying reduced subscheme.

\begin{defi} 
The intersection graph of $E_{\rm red}$ is defined as 
the following undirected graph:\smallskip

-vertices. A vertex of multiplicity $n$ for each irreducible 
component $(E_{\rm red})_i$ of $E_{\rm red}$ which decomposes 
over the algebraic closure of $K$ into $n$ irreducible 
components.\smallskip

-edges. Different $(E_{\rm red})_i$ and $(E_{\rm red})_j$ are 
connected by $(E_{\rm red})_i.(E_{\rm red})_j$ 
undirected edges. \\ $(E_{\rm red})_i$ has\ \ 
$\frac{1}{2}(E_{\rm red})_i.(E_{\rm red})_i+
\textit{multiplicity of $(E_{\rm red})_i$}$\ \ loops.
\end{defi}

If $K$ is algebraically closed, then the $(E_{\rm red})_i$ are
isomorphic to $\P^1_K$ and the self-intersection of each 
$(E_{\rm red})_i$ is $-2$, because the resolution is crepant.\medskip

Let $K\to L$ be a Galois extension, $\Gamma=\Aut_K(L)$.
$\Gamma$ operates on the intersection graph of $(E_{\rm red})_L$ 
by graph automorphisms. The irreducible components 
$(E_{\rm red})_i$ of $E_{\rm red}$ correspond to $\Gamma$-orbits 
of irreducible components $(E_{\rm red})_{L,k}$ of 
$(E_{\rm red})_L$ by proposition \ref{prop:irredcomp}. 
For the intersection form one has 
\[\textstyle (E_{\rm red})_i.(E_{\rm red})_j=
((E_{\rm red})_i)_L.((E_{\rm red})_j)_L=
\sum_{kl}(E_{\rm red})_{L,k}.(E_{\rm red})_{L,l}\]
where indices $k$ and $l$ run through the irreducible 
components of $((E_{\rm red})_i)_L$ and $((E_{\rm red})_j)_L$ 
respectively. Thus for the intersection graph there is a 
proposition similar to proposition \ref{prop:repgraph} for 
representation graphs.  

\begin{prop} \label{prop:intgraph}
The intersection graph of $E_{\rm red}$ arises by identifying 
the elements of $\Gamma$-orbits of vertices of the intersection 
graph of $(E_{\rm red})_L$, adding multiplicities. 
The edges between vertices $(E_{\rm red})_i$ and $(E_{\rm red})_j$
are in bijection with the edges between the irreducible 
components of $((E_{\rm red})_i)_L$ and $((E_{\rm red})_j)_L$.\qed
\end{prop}

\subsection{Irreducible components of $E$ and irreducible 
representations of $G$}

The basic statement of McKay correspondence is a bijection
between the set of irreducible components 
of the exceptional divisor $E$ and the set of isomorphism classes 
of nontrivial irreducible representations of the group scheme $G$.

\begin{thm}\label{thm:McKaybij} 
There are bijections for intermediate fields $K\sub L\sub C$
between the set $\Irr(E_L)$ of irreducible components of $E_L$ 
and the set $\Irr(G_L)$ of isomorphism classes of nontrivial 
irreducible representations of $G_L$ having the property
that for $K\sub L\sub L'\sub C$, if the bijection 
$\Irr(E_L)\to\Irr(G_L)$ for $L$ maps $E_i\mapsto V_i$, then 
the bijection $\Irr(E_{L'})\to\Irr(G_{L'})$ for $L'$ maps 
irreducible components of $(E_i)_{L'}$ to irreducible 
components of $(V_i)_{L'}$.
\end{thm}

\begin{proof} As described earlier, the Galois group 
$\Gamma=\Aut_L(C)$ of the Galois extension $L\to C$, operates 
on the sets $\Irr(G_C)$ and $\Irr(E_C)$. In both cases elements 
of $\Irr(G_L)$ and $\Irr(E_L)$ correspond to $\Gamma$-orbits of 
elements of $\Irr(G_C)$ and $\Irr(E_C)$ by corollary 
\ref{cor:irredrep} and proposition \ref{prop:irredcomp} 
respectively.
This way a given bijection between the sets $\Irr(G_C)$ and 
$\Irr(E_C)$ defines a bijection between $\Irr(G_L)$ and 
$\Irr(E_L)$ on condition that the bijection is equivariant with 
respect to the operations of $\Gamma$. 
Checking this for the bijection of McKay correspondence
over the algebraically closed field $C$ constructed via
stratification or via tautological sheaves will give bijections
over intermediate fields $L$ having the property of the theorem.
This will be done in the process of proving theorem 
\ref{thm:strati} or theorem \ref{thm:tautol}. 
\end{proof}

Moreover, in the situation of the theorem the Galois group
$\Gamma=\Aut_L(C)$ operates on the representation graph of 
$G_C$ and on the intersection graph of $(E_{\rm red})_C$. 
Then in both cases the graphs over $L$ arise by identifying 
the elements of $\Gamma$-orbits of vertices of the graphs 
over $C$ by proposition \ref{prop:repgraph} and 
\ref{prop:intgraph}. 
Therefore an isomorphism of the graphs over $C$, the bijection 
between the sets of vertices being $\Gamma$-equivariant,
defines an isomorphism of the graphs over $L$.

For the algebraically closed field $C$ this is the classical 
McKay correspondence for subgroups of $\SL(2,C)$ (\cite{McK80},
\cite{GV83}, \cite{ItNm99}).
The statement, that there is a bijection of edges between
given vertices $(E_{\rm red})_{L,i}\leftrightarrow V_i$ and
$(E_{\rm red})_{L,j}\leftrightarrow V_j$, can be formulated
equivalently in terms of the intersection form as 
$(E_{\rm red})_{L,i}.(E_{\rm red})_{L,j}=\<V_i,V_j\>$.

\begin{thm}\label{thm:McKaygraph}
The bijections $E_i\leftrightarrow V_i$ of theorem \ref{thm:McKaybij}
between irreducible components of $E_L$ and isomorphism classes 
of nontrivial irreducible representations of $G_L$ can be 
constructed such that $(E_{\rm red})_i.(E_{\rm red})_j=\<V_i,V_j\>$
or equivalently that these bijections define isomorphisms of 
graphs between the intersection graph of $(E_{\rm red})_L$
and the representation graph of $G_L$.\qed
\end{thm}

We will consider two ways to construct bijections between
nontrivial irreducible representations and irreducible
components with the properties of theorem \ref{thm:McKaybij} 
and \ref{thm:McKaygraph}: A stratification of $\GHilb_K\A^2_K$ 
(\cite{ItNm96}, \cite{ItNm99})
and the tautological sheaves on $\GHilb_K\A^2_K$
(\cite{GV83}, \cite{KaVa00}).

\subsection{Stratification of $\GHilb_K\A^2_K$}

Let $S:=K[x_1,x_2]$, let $O\in\A_K^2$ be the origin, 
$\mathfrak m\subset S$ the corresponding maximal ideal, 
$\overline{O}:=\pi(O)\in\A_K^2/G$ with corresponding
maximal ideal $\mathfrak n\subset S^G$, let $\overline{S}:=
S/\mathfrak n S$ with maximal ideal $\overline{\mathfrak m}$.
An $L$-valued point of the fiber $E=\tau^{-1}(\overline{O})$ 
corresponds to a $G$-cluster defined by an ideal $I\subset S_L$ 
such that $\mathfrak n_L\sub I$ or equivalently an ideal
$\overline{I}\subset\overline{S}_L=S_L/\mathfrak n_L S_L$.
For such an ideal $I$ define the 
representation $V(I)$ over $L$ by
\[V(I):=\overline{I}/\overline{\mathfrak m}_L\overline{I}\]

\begin{lemma}\label{le:gammaI} 
For $\gamma\in\Aut_K(L)$:\ $V(\gamma^{-1}I)\cong V(I)^\gamma$.
\end{lemma}
\begin{proof}
As an $A_L$-comodule $\overline{I}=\overline{I}_0\oplus
\overline{\mathfrak m}_L\overline{I}$, where 
$\overline{I}_0\cong \overline{I}/\overline{\mathfrak m}_L\overline{I}$.
Then $\gamma^{-1}\overline{I}=\gamma^{-1}\overline{I}_0\oplus
\overline{\mathfrak m}_L(\gamma^{-1}\overline{I})$
and $V(\gamma^{-1}I)=\gamma^{-1}\overline{I}/\overline{\mathfrak m}_L
(\gamma^{-1}\overline{I})\cong\gamma^{-1}\overline{I}_0
\cong\overline{I}_0^\gamma\cong V(I)^\gamma$ by remark 
\ref{rem:gammaU} applied to $\overline{I}_0\sub\overline{S}_L$.
\end{proof}

\begin{thm}\label{thm:strati} 
There is a bijection $E_j\leftrightarrow V_j$ between the set 
$\Irr(E)$ of irreducible components of $E$ and the set 
$\Irr(G)$ of isomorphism classes of nontrivial irreducible 
representations of $G$ such that for any closed point $y\in E$:
If $I\subset S_{\kappa(y)}$ is an ideal defining a $\kappa(y)$-valued 
point of the scheme $\{y\}\subset E$, then 
\[
\Hom^G_{\kappa(y)}(V(I),(V_j)_{\kappa(y)})\neq 0\;
\Longleftrightarrow\; y\in E_j
\]
and $V(I)$ is either irreducible or consists of two irreducible 
representations not isomorphic to each other.
Applied to the situation after base extension $K\to L$, $L$ an 
algebraic extension of $K$, one obtains bijections 
$\Irr(E_L)\leftrightarrow\Irr(G_L)$ having the
properties of theorems \ref{thm:McKaybij} and \ref{thm:McKaygraph}.
\end{thm}
\begin{proof}
In the case of algebraically closed $K$ the theorem follows from
\cite{ItNm99}.

In the general case denote by $U_i$ the isomorphism classes 
of nontrivial irreducible representations of $G_C$ over 
the algebraic closure $C$.
Over $C$ the theorem is valid, let $E_{C,i}$ be the component 
corresponding to $U_i$.\\
We show that this bijection is equivariant with respect to the
operations of $\Gamma=\Aut_K(C)$.
Let $x\in E_{C,i}$ be a closed point such that $x\not\in E_{C,i'}$ 
for $i'\neq i$.
Then for the corresponding $C$-valued point 
$\alpha\colon\Spec C\to E_{C,i}$ given by an ideal $I\subset S_C$ one 
has $V(I)\cong U_i$. 
By corollary \ref{cor:conjpoint} the $C$-valued point corresponding
to $\gamma x$ is $\alpha^\gamma$ given by the ideal
$\gamma^{-1}I\subset S_C$.
By lemma \ref{le:gammaI} $V(\gamma^{-1}I)\cong U_{\gamma(i)}$, 
where $U_{\gamma(i)}=U_i^\gamma$. 
Therefore $\gamma x\in E_{\gamma(i)}$ and $\gamma E_i=E_{\gamma(i)}$.\\
For an irrreducible representation $V_j$ of $G$ define 
$E_j$ to be the component of $E$, which decomposes over $C$
into the irreducible components $E_{C,i}$ satisfying 
$U_i\sub (V_j)_C$. This method, applied to the situation after 
base extension $K\to L$, leads to bijections having the properties 
of theorems \ref{thm:McKaybij} and \ref{thm:McKaygraph}.

We show that this bijection is given by the condition in the 
theorem. 
Let $y$ be a closed point of $E$ and $\alpha$ a $\kappa(y)$-valued 
point of the scheme $\{y\}\subset E$ given by an ideal 
$I\subset S_{\kappa(y)}$. 
$K\to\kappa(y)$ is an algebraic extension, embed $\kappa(y)$ into $C$.
After base extension $\kappa(y)\to C$ one has the $C$-valued point
$\alpha_C\colon\Spec C\to\{y\}_C$ given by $I_C\subset S_C$.
Then $V(I)_C\cong V(I_C)$ and $I_C$ corresponds to a closed point 
$z\in\{y\}_C\subset E_C$. Therefore
\[\begin{array}{rcl}
y\in E_j&\Longleftrightarrow& z\in E_{C,i}\;\;
\textit{for some $i$ satisfying $U_i\sub (V_j)_C$}\\
&\Longleftrightarrow&\Hom^G_C(V(I_C),U_i)\neq 0\;\;
\textit{for some $i$ satisfying $U_i\sub (V_j)_C$}\\
&\Longleftrightarrow&\Hom^G_{\kappa(y)}(V(I),(V_j)_{\kappa(y)})\neq 0\\
\end{array}
\]\vspace{-2.2em}

\end{proof}

\subsection{Tautological sheaves}

Let $0\to\mathscr I\to\O_{\A^2_Y}\to\O_Z\to 0$
be the universal quotient of $Y:=\GHilb_K\A^2_K$.
The projection $p\colon Z\to Y$ is a finite flat morphism,
$p_*\O_Z$ is a locally free $G$-sheaf on $Y$
with fibers $p_*\O_Z\otimes_{\O_Y}\kappa(y)$ isomorphic to the
regular representation over $\kappa(y)$.

Let $V_0,\ldots,V_s$ the isomorphism classes of irreducible 
representations of $G$, $V_0$ the trivial representation.
The $G$-sheaf $\mathscr G:=p_*\O_Z$ on $Y$ decomposes into 
isotypic components (see remark \ref{rem:Gsheaves}.(3) and 
subsection \ref{subsec:decomp-Galoisext})
\[\textstyle\mathscr G\cong\bigoplus_{j=0}^s\mathscr G_j\]
where $\mathscr G_j$ is the component for $V_j$.

\begin{defi}
For any isomorphism class $V_j$ of irreducible representations
of $G$ over $K$ define the sheaf $\mathscr F_j$ on 
$Y=\GHilb_K\A^2_K$ by
\[
\mathscr F_j:=\SHom^G_{\O_Y}(V_j\otimes_K\O_Y,\mathscr G_j)
=\SHom^G_{\O_Y}(V_j\otimes_K\O_Y,\mathscr G)
\]
For a field extension $K\to L$ denote by $\mathscr F_{L,i}$ 
the sheaf $\SHom^{G_L}_{\O_{Y_L}}(U_i\otimes_L\O_{Y_L},\mathscr G_L)$
on $Y_L$, $U_i$ an irreducible representation of $G_L$ over $L$.
\end{defi}

\begin{rem}\ \\
(1) For $K=\C$ the sheaves $\mathscr F_j$ were studied in 
\cite{GV83}, \cite{KaVa00}, they may be defined as well as
$\mathscr F_j=\tau^*\SHom^G_{\A^2_K/G}(V_j\otimes_K\O_{\A^2_K/G},
\pi_*\O_{\A^2_K})/(\O_Y\textit{-torsion})$ or 
$(p_*q^*(\O_{\A^2_K}\otimes_K V_j^\vee))^G$ 
using the canonical morphisms in the diagram

\vspace{-4.5mm}

\begin{picture}(160,30)
\put(80,25){\makebox(0,0)[c]{$Z$}}
\put(60,15){\makebox(0,0)[c]{$Y$}}
\put(100,15){\makebox(0,0)[c]{$\A^2_K$}}
\put(80,5){\makebox(0,0)[c]{$\A^2_K/G$}}
\put(77,23){\vector(-2,-1){13}}
\put(83,23){\vector(2,-1){13}}
\put(63,13){\vector(2,-1){11}}
\put(97,13){\vector(-2,-1){11}}
\put(69,21.5){\makebox(0,0)[c]{\footnotesize$p$}}
\put(91,21.5){\makebox(0,0)[c]{\footnotesize$q$}}
\put(68,8.5){\makebox(0,0)[c]{\footnotesize$\tau$}}
\put(92,8.5){\makebox(0,0)[c]{\footnotesize$\pi$}}
\end{picture}

\vspace{-1.5mm}

(2) $\mathscr F_j$ is a locally free sheaf of rank $\dim_K V_j$.\\
(3) For each $j$ there is the natural isomorphism of $G$-sheaves
$\mathscr F_j\otimes_{\End^G_K(V_j)}V_j
\stackrel{\!_\sim}{\longrightarrow}\mathscr G_j$.
\end{rem}

Let $K\to L$ be a Galois extension and $U_0,\ldots,U_r$ be the 
isomorphism classes of irreducible representations of $G_L$ over $L$. 
Then a decomposition $(V_j)_L=\bigoplus_{i\in I_j}U_i$ over $L$ 
of an irreducible representation $V_j$ of $G$ over $K$ gives 
a decomposition of the corresponding tautological sheaf
\[
\begin{array}{l}
(\mathscr F_j)_L=\SHom^G_{\O_Y}(V_j\otimes_K\O_Y,\mathscr G)_L
\cong\SHom^{G_L}_{\O_{Y_L}}((V_j\otimes_K\O_Y)_L,\mathscr G_L)\\
\cong\SHom^{G_L}_{\O_{Y_L}}(\bigoplus_{i\in I_j}U_i\otimes_L\O_{Y_L},
\mathscr G_L)\cong\bigoplus_{i\in I_j}\SHom^{G_L}_{\O_{Y_L}}
(U_i\otimes_L\O_{Y_L},\mathscr G_L)
=\bigoplus_{i\in I_j}\mathscr F_{L,i}\\
\end{array}
\]
We have used the fact that the $U_i$ occur with multiplicity $1$ 
as it is the case for finite subgroup schemes of $\SL(2,K)$, see 
proposition \ref{prop:m}.\smallskip

The tautological sheaves $\mathscr F_j$ can be used to establish 
a bijection between the set of irreducible components of 
$E_{\rm red}$ and the set of isomorphism classes of nontrivial
irreducible representations of $G$ by considering intersections 
$\mathscr L_j.(E_{\rm red})_{j'}$, i.e. the degrees of 
restrictions of the line bundles
$\mathscr L_j:=\bigwedge^{\rk\mathscr F_j}\mathscr F_j$ 
to the curves $(E_{\rm red})_{j'}$.

\begin{thm}\label{thm:tautol}
There is a bijection $E_j\leftrightarrow V_j$ between the set 
$\Irr(E)$ of irreducible components of $E$ and the set 
$\Irr(G)$ of isomorphism classes of nontrivial irreducible 
representations of $G$ such that

\vspace{-4.5mm}

\[\mathscr L_j.(E_{\rm red})_{j'}=\dim_K\Hom^G_K(V_j,V_{j'})\]

\vspace{-2mm}

where $\mathscr L_j=\bigwedge^{\rk\mathscr F_j}\mathscr F_j$.\\
Applied to the situation after base extension $K\to L$, $L$ an 
algebraic extension field of $K$, one obtains bijections 
$\Irr(E_L)\leftrightarrow\Irr(G_L)$ having the
properties of theorems \ref{thm:McKaybij} and \ref{thm:McKaygraph}.
\end{thm}
\begin{proof}
In the case of algebraically closed $K$ the theorem follows from
\cite{GV83}.

In the general case denote by $U_0,\ldots, U_r$ the isomorphism 
classes of irreducible representations of $G_C$ over the algebraic 
closure $C$, $U_0$ the trivial one.
Over $C$ the theorem is  valid, let $E_{C,i}$ be the component 
corresponding to $U_{i}$, what means that 
$\mathscr L_{C,i}.(E_{\rm red})_{C,i'}=\delta_{ii'}$,
where $\mathscr L_{C,i}=\bigwedge^{\rk\mathscr F_{C,i}}
\mathscr F_{C,i}$.

To show that the bijection over $C$ is equivariant with respect 
to the operations of $\Gamma=\Aut_K(C)$, one has to show that
$\gamma_*\mathscr L_{C,i}\cong\mathscr L_{C,\gamma(i)}$,
where $U_{\gamma(i)}=U_i^\gamma$.
Then $\mathscr L_{C,i}.E_{C,i'}
=\gamma_*\mathscr L_{C,i}.\gamma E_{C,i'}
=\mathscr L_{C,\gamma(i)}.\gamma E_{C,i'}$
and therefore $\gamma E_{C,i'}=E_{C,\gamma(i')}$.
It is $\gamma_*\mathscr L_{C,i}\cong\mathscr L_{C,\gamma(i)}$,
because using lemma \ref{le:gamma_*} and remark \ref{rem:F->gamma_*F}

\vspace{-3mm}

\[
\gamma_*\mathscr F_{C,i}
\cong\SHom^{G_C}_{\O_{Y_C}}(\gamma_*(U_i\otimes_C\O_{Y_C}),
\gamma_*\mathscr G_C)\cong\SHom^{G_C}_{\O_{Y_C}}
(U_i^\gamma\otimes_C\O_{Y_C},\mathscr G_C)=\mathscr F_{C,\gamma(i)}
\]

\vspace{-2mm}

Since the bijection over $C$ is equivariant with respect to the 
$\Gamma$-operations on $\Irr(G_C)$ and $\Irr(E_C)$, one
can define a bijection $\Irr(G)\leftrightarrow\Irr(E)$:
For $V_j\in\Irr(G)$ let $E_j$ be the element of $\Irr(E)$
such that $(V_j)_C=\bigoplus_{i\in I_j}U_i$ and
$(E_j)_C=\bigcup_{i\in I_j}E_{C,i}$ for the same subset 
$I_j\sub\{1,\ldots,r\}$. This method applied to the situation after
base extension $K\to L$ leads to bijections having the properties 
of theorems \ref{thm:McKaybij} and \ref{thm:McKaygraph}.

We show that this bijection is given by the construction of 
the theorem. It is $(\mathscr F_j)_C=\bigoplus_{i\in I_j}
\mathscr F_{C,i}$ and therefore

\vspace{-5.5mm}

\[
\begin{array}{l}
\mathscr L_j.(E_{\rm red})_{j'}=(\mathscr L_j)_C.((E_{\rm red})_{j'})_C
=\big(\bigotimes_{i\in I_j}\mathscr L_{C,i}\big).
\big(\sum_{i'\in I_{j'}}(E_{\rm red})_{C,i'}\big)
=\sum_{i,i'}\mathscr L_{C,i}.(E_{\rm red})_{C,i'}\\
=\sum_{i,i'}\dim_C\Hom^{G_C}_C(U_i,U_{i'})
=\dim_C\Hom^{G_C}_C((V_j)_C,(V_{j'})_C)
=\dim_K\Hom^G_K(V_j,V_{j'})\\
\end{array}
\]

\vspace{-6mm}
\end{proof}

\subsection{Examples}
\label{subsec:ex}

{\bf Finite subgroups of \boldmath $\SL(2,K)$.}
In the case of subgroups $G\subset\SL(2,K)$ the representation
graph can be read off from the table of characters of the
group $G$ over an algebraically closed field, since in this case
representations are conjugate if and only if the values of 
their characters are.
We have the following graphs for the finite subgroups of 
$\SL(2,K)$:\\
- cyclic group $\Z/n\Z$, $n\geq 1$.
It is $\xi+\xi^{-1}\in K$, $\xi$ a primitive $n$-th root of unity.
Diagram $(A_{n-1})$ if $\xi\in K$, otherwise $(A_{n-1})'$.\\ 
- binary dihedral group $\BD_n$, $n\geq 2$. 
It is $c=\frac{1}{2}(\xi+\xi^{-1})\in K$, $\xi$ a primitive $2n$-th 
root of unity, and $\HSymb{-1}{c^2-1}_K=1$.
Diagram $(D_{n+2})$ if $n$ even or $\sqrt{-1}\in K$, otherwise 
$(D_{n+2})'$.\\
- binary tetrahedral group $\BT$. It is $\HSymb{-1}{-1}_K=1$.
Diagram $(E_6)$ if $K$ contains a primitive $3$rd root of unity,
otherwise $(E_6)'$.\\ 
- binary octahedral group $\BO$. 
It is $\HSymb{-1}{-1}_K=1$ and $\sqrt{2}\in K$.
Diagram $(E_7)$.\\
- binary icosahedral group $\BI$. 
It is $\HSymb{-1}{-1}_K=1$ and $\sqrt{5}\in K$.
Diagram $(E_8)$.\smallskip

Examples for the graphs $(A_n)'$, $(D_{2m+1})'$, $(E_6)'$:\\
\begin{tabular}{ll}
$(A_n)'$& $\Z/(n+1)\Z$ over $\Q(\xi+\xi^{-1})$, 
$\xi$ a primitive $(n+1)$-th root of unity.\\
$(D_{2m+1})'$& $\BD_{2m-1}$ over $\Q(\xi)$, 
$\xi$ a primitve $2(2m-1)$-th root of unity.\\
$(E_6)'$& $\BT$ over $\Q(\sqrt{-1})$.\\
\end{tabular}

\medskip

{\bf Abelian subgroup schemes.}
In the case of abelian subgroup schemes of $\SL(2,K)$ the
graphs $(A_n)$ and $(A_n)'$ occur.\\
- the cyclic group $G=\Z/n\Z$ is realisable as the subgroup of
$\SL(2,K)$ generated by $g:=\left(\begin{smallmatrix}0&-1\\
1&\xi+\xi^{-1}\end{smallmatrix}\right)$, if the field $K$ 
contains $\xi+\xi^{-1}$ for $\xi$ a primitive $n$-th root of 
unity. 
If $K$ does not contain $\xi$, then there are $1$-dimensional 
representations over the algebraic closure that are not
realisable over $K$, one has diagram $(A_{n-1})'$.\\
- for the subgroup scheme $G=\mu_n\subset\SL(2,K)$ the Hopf 
algebra $K[y]/\<y^n\>$ decomposes into a direct sum of 
simple subcoalgebras $\<y^j\>_K$ corresponding to $1$-dimensional
representations of $G$. Thus one has diagram $(A_{n-1})$.

\medskip

{\bf The graph \boldmath$(D_{2m})'$.}
Let $n\geq 2$, $\e$ a primitive $4n$-th root of unity and
$\xi=\e^2$. Put $K=\Q(\e+\e^{-1})$, $C=\Q(\e)$ and 
$\Gamma=\Aut_K(C)=\{\id,\gamma\}$. One has the injective 
representation of 
$\BD_n=\<\sigma,\tau\:|\:\tau^2=\sigma^n=(\tau\sigma)^2\>$ 
in $\SL(2,C)$ given by

\vspace{-2.5mm}

\[
\sigma\mapsto\left(\begin{matrix}\xi&0\\0&\xi^{-1}\end{matrix}\right),
\qquad
\tau\mapsto\left(\begin{matrix}0&-\e\\\e^{-1}&0\end{matrix}\right)
\]

\vspace{-1mm}

We will identify $\BD_n$ with its image in $\SL(2,C)$ and regard 
it as a subgroup scheme of $\SL(2,C)$.\smallskip

$\Gamma$ operates on $\SL(2,C)$, the $K$-automorphism
$\gamma\in\Gamma$, $\gamma\colon\e\mapsto\e^{-1}$ of order $2$ 
operates nontrivially on the closed points of $\BD_n$ by
$\sigma\mapsto\sigma^{-1}$, $\tau\mapsto\tau\sigma$.
The subgroup scheme $\BD_n\subset\SL(2,C)$ is defined over $K$,
let $G\subset\SL(2,K)$ such that $G_C=\BD_n$. 
The closed points of $G$ correspond to $\Gamma$-orbits of closed 
points of $\BD_n$, they have the form 
$\{\id\}$, $\{-\id\}$, $\{\sigma^k,\sigma^{-k}\}$, 
$\{\tau\sigma^k,\tau\sigma^{-k+1}\}$.\smallskip

\enlargethispage{3mm}

The automorphism $\gamma$ operates on the characters of $\BD_n$ 
trivially except that for even $n$ it permutes two of the 
irreducible $1$-dimensional representations.
One has the graph $(D_{n+2})'$ for $n$ even and the graph $(D_{n+2})$ 
for $n$ odd.

\newpage
\appendix
\section{Finite subgroups of $\SL(2,C)$: Presentations and character
tables}\label{app:charsubgrSL2}

- Cyclic groups\medskip

The irreducible representations are
$\chi_j\colon\Z/n\Z\to C^*$, $\overline{i}\mapsto\xi^{ji}$
for $j\in\{0,\ldots,n-1\}$,\\ 
where $\xi$ is a primitive $n$-th root of unity.\bigskip\medskip

- Binary dihedral groups:\ \
$BD_n=\<\sigma,\tau\:|\:\tau^2=\sigma^n=(\tau\sigma)^2\>$,\ \ 
$-\id:=(\tau\sigma)^2$.\smallskip

{\footnotesize
\begin{tabular}{lcl}
$BD_n$, $n$ odd&&$BD_n$, $n$ even\\
\begin{tabular}{l|ccccc}
&$\id$&$-\id$&$\sigma^k$&$\tau$&$\tau\sigma$\\
\hline
$\mathbf{1}$&$1$&$1$&$1$&$1$&$1$\\
$\mathbf{1}'$&$1$&$1$&$1$&$-1$&$-1$\\
$\mathbf{1}''$&$1$&$-1$&$(-1)^k$&$i$&$-i$\\
$\mathbf{1}'''$&$1$&$-1$&$(-1)^k$&$-i$&$i$\\
$\mathbf{2}^j$&$2$&$(-1)^j2$&$\xi^{kj}+\xi^{-kj}$&$0$&$0$\\
\end{tabular}
&$\quad$&
\begin{tabular}{l|ccccc}
&$\id$&$-\id$&$\sigma^k$&$\tau$&$\tau\sigma$\\
\hline
$\mathbf{1}$&$1$&$1$&$1$&$1$&$1$\\
$\mathbf{1}'$&$1$&$1$&$1$&$-1$&$-1$\\
$\mathbf{1}''$&$1$&$1$&$(-1)^k$&$1$&$-1$\\
$\mathbf{1}'''$&$1$&$1$&$(-1)^k$&$-1$&$1$\\
$\mathbf{2}^j$&$2$&$(-1)^j2$&$\xi^{kj}+\xi^{-kj}$&$0$&$0$\\
\end{tabular}\\
\end{tabular}\smallskip

$\xi$ a primitive $2n$-th root of unity and $j=1,\ldots,n-1$.}\bigskip\medskip

- Binary tetrahedral group:\ \
$BT=\<a,b\:|\:a^3=b^3=(ab)^2\>$,\ \ $-id:=(ab)^2$.\smallskip

{\footnotesize
\begin{tabular}{l|ccccccc|c}
&$\id$&$-\id$&$a$&$-a$&$b$&$-b$&$ab$\\
\hline
${\mathbf 1}$&$1$&$1$&$1$&$1$&$1$&$1$&$1$&$1$\\
${\mathbf 1'}$&$1$&$1$&$\omega$&$\omega$&$\omega^2$&$\omega^2$&$1$&$\Z/3\Z$\\
${\mathbf 1''}$&$1$&$1$&$\omega^2$&$\omega^2$&$\omega$&$\omega$&$1$&$\Z/3\Z$\\
${\mathbf 3}$&$3$&$3$&$0$&$0$&$0$&$0$&$-1$&$A_4$\\
${\mathbf 2}$&$2$&$-2$&$1$&$-1$&$1$&$-1$&$0$&$BT$\\
${\mathbf 2'}$&$2$&$-2$&$\omega$&$-\omega$&$\omega^2$&$-\omega^2$&$0$&$BT$\\
${\mathbf 2''}$&$2$&$-2$&$\omega^2$&$-\omega^2$&$\omega$&$-\omega$&$0$&$BT$\\
\hline
&$1$&$1$&$4$&$4$&$4$&$4$&$6$&
\end{tabular}\\
$\omega$ a primitive $3$rd root of unity.}\bigskip\medskip

- Binary octahedral group:\ \
$BO=\<a,b\:|\:a^3=b^4=(ab)^2\>$,\ \ $-\id:=(ab)^2$.\smallskip

{\footnotesize
\begin{tabular}{l|cccccccc|c}
&$\id$&$-\id$&$ab$&$a$&$-a$&$b$&$-b$&$b^2$&\\
\hline
${\mathbf 1}$&$1$&$1$&$1$&$1$&$1$&$1$&$1$&$1$&$1$\\
${\mathbf 1'}$&$1$&$1$&$-1$&$1$&$1$&$-1$&$-1$&$1$&$\Z/2\Z$\\
${\mathbf {2'''}}$&$2$&$2$&$0$&$-1$&$-1$&$0$&$0$&$2$&$S_3$\\
${\mathbf 3}$&$3$&$3$&$1$&$0$&$0$&$-1$&$-1$&$-1$&$S_4$\\
${\mathbf 3'}$&$3$&$3$&$-1$&$0$&$0$&$1$&$1$&$-1$&$S_4$\\
${\mathbf 2}$&$2$&$-2$&$0$&$1$&$-1$&$\sqrt{2}$&$-\sqrt{2}$&$0$&$BO$\\
${\mathbf 2'}$&$2$&$-2$&$0$&$1$&$-1$&$-\sqrt{2}$&$\sqrt{2}$&$0$&$BO$\\
${\mathbf 4}$&$4$&$-4$&$0$&$-1$&$1$&$0$&$0$&$0$&$BO$\\
\hline
&$1$&$1$&$12$&$8$&$8$&$6$&$6$&$6$&\\
\end{tabular}}\bigskip\medskip

- Binary icosahedral group:\ \
$BI=\<a,b\:|\:a^3=b^5=(ab)^2\>$,\ \ $-\id:=(ab)^2$.\smallskip

{\footnotesize
\begin{tabular}{l|ccccccccc|c}
&$\id$&$-\id$&$a$&$-a$&$b$&$-b$&$b^2$&$-b^2$&$ab$\\
\hline
$\mathbf 1$&$1$&$1$&$1$&$1$&$1$&$1$&$1$&$1$&$1$&$1$\\
$\mathbf 3$&$3$&$3$&$0$&$0$&$\mu^+$&$\mu^+$&$\mu^-$&$\mu^-$&$-1$&$A_5$\\
$\mathbf 3'$&$3$&$3$&$0$&$0$&$\mu^-$&$\mu^-$&$\mu^+$&$\mu^+$&$-1$&$A_5$\\
$\mathbf 4'$&$4$&$4$&$1$&$1$&$-1$&$-1$&$-1$&$-1$&$0$&$A_5$\\
$\mathbf 5$&$5$&$5$&$-1$&$-1$&$0$&$0$&$0$&$0$&$0$&$A_5$\\

$\mathbf 2$&$2$&$-2$&$1$&$-1$&$\mu^+$&$-\mu^+$&$-\mu^-$&$\mu^-$&$0$&$BI$\\
$\mathbf 2'$&$2$&$-2$&$1$&$-1$&$\mu^-$&$-\mu^-$&$-\mu^+$&$\mu^+$&$0$&$BI$\\
$\mathbf 4$&$4$&$-4$&$-1$&$1$&$1$&$-1$&$-1$&$1$&$0$&$BI$\\
$\mathbf 6$&$6$&$-6$&$0$&$0$&$-1$&$1$&$1$&$-1$&$0$&$BI$\\

\hline
&$1$&$1$&$20$&$20$&$12$&$12$&$12$&$12$&$30$\\
\end{tabular}\\
$\mu^+:=\frac{1}{2}(1+\sqrt{5})$, $\mu^-:=\frac{1}{2}(1-\sqrt{5})$.}

\bigskip\bigskip

Mark Blume\\
Mathematisches Institut, Universit\"at M\"unster\\
Einsteinstrasse 62, 48149 M\"unster, Germany\\
Email: mark.blume@uni-muenster.de

\end{document}